\documentclass[twocolumn]{autart}

\usepackage{amssymb}
\usepackage{amsmath}
\usepackage{graphicx}
\graphicspath{ {fig/} }
\usepackage{epstopdf}
\usepackage{csquotes}

\begin{document}

\begin{frontmatter}

\title{Consensus and Formation Control on $SE(3)$ for Switching Topologies \thanksref{footnoteinfo}}

\thanks[footnoteinfo]{The authors gratefully acknowledge the financial support form the Fonds
National de la Recherche, Luxembourg (8864515)}

\author[Johan]{Johan Thunberg}\ead{johan.thunberg@uni.lu},
\author[Johan]{Jorge Goncalves}\ead{jorge.goncalves@uni.lu},
\author[Hu]{Xiaoming Hu}\ead{hu@math.kth.se}

\address[Johan]{Universite du Luxembourg
Campus Belval
7, avenue des Hauts-Fourneaux
L-4362 Esch-sur-Alzette} 

\address[Hu]{KTH Royal Institute of Technology,  S-100 44 Stockholm, Sweden}


\begin{abstract}
This paper addresses the consensus problem and the formation problem on $SE(3)$ in multi-agent systems
with directed and switching interconnection topologies.
Several control laws are introduced for the consensus problem. By a simple transformation, it is shown that the proposed control laws can be used
for the formation problem. The design is first conducted on the kinematic level, where the velocities are the control laws. Then, for rigid bodies in space,
the design is conducted on the dynamic level, where the torques and the forces are the control laws. On the kinematic level, first two 
control laws are introduced that explicitly use
Euclidean transformations, then 
separate control laws are defined for the rotations and the translations.
In the special case of purely rotational motion, the consensus problem
is referred to as \emph{consensus on $SO(3)$} or \emph{attitude synchronization}. In this problem,
for a broad class of
local representations or parameterizations of $SO(3)$, including the Axis-Angle Representation, the
Rodrigues Parameters and the Modified Rodrigues Parameters, two types of control laws
are presented that look structurally the same for any choice of local representation. For these two control laws we provide conditions on the initial rotations and the connectivity of the graph such that the system reaches consensus on $SO(3)$. Among the contributions of this 
paper, there are conditions for when exponential rate of convergence occur. A theorem is provided showing that for any choice of local representation for the rotations, there is a change of coordinates such that the transformed system has a well known structure.
\end{abstract}

\end{frontmatter}

\section{Introduction}
This work addresses the problem of continuous time consensus and formation control on $SE(3)$ for switching interconnection topologies.
We start by designing control laws
on a kinematic level, where the velocities are
control signals and then continue to design
control laws on the dynamic level for rigid bodies in space,
 where the forces and 
the torques are the control laws. 

The main 
focus in this paper is on the kinematic control laws. 
This approach is not justified from a physical perspective. Nevertheless, there are reasons why this 
path is still reasonable to take. 
Firstly, dynamics are often platform dependent, and
especially in the robotics community it is desired
to specify control laws on the kinematic level. Secondly,
a deeper understanding of how the geometry of $SE(3)$ (and in particular $SO(3)$) affects the 
control design can be acquired by designing the control 
laws on a kinematic level, since we are then 
working directly in the tangent space of $SE(3)$ or
$SO(3)$. 

On $SE(3)$, consensus control is a special case of formation control, but actually formation control 
can be seen as a special case of consensus control,
a fact that will be used in this paper. The approach 
is to develop consensus control laws,
which after a simple transformation can be used as 
formation control laws. 
By taking this approach we can use existing theory 
for consensus in order to provide convergence results
for the formation problem.

The consensus control problem on $SO(3)$ comprises a subset
of the consensus control problem on $SE(3)$, but it
is, from many perspectives, the most challenging part of the control design. Hence, most emphasis will be taken towards this problem. 
Whereas the translations are elements of $\mathbb{R}^3$, the rotations are elements of the compact manifold $SO(3)$, the group of orthogonal matrices in $\mathbb{R}^{3\times3}$ with determinant equal to~$1$. 

There is a wide range of applications for the proposed  
control laws, \emph{e.g.}, satellites or spacecraft  that shall reach a certain formation, multiple robotic arms that shall hold a rigid object or
cameras that shall look in some desired 
directions (or the same direction in case of consensus). 
For rigid bodies in space, \emph{i.e.}, spacecraft or satellites,
there has recently been an extensive research on the consensus on $SO(3)$ problem~\cite{listmann2009passivity,dimarogonas2006laplacian,sarlette2009automatica,ren2010distributed,tron2012intrinsic}. In that problem, the goal is
to design a control torque such that the rotations
of the rigid bodies become synchronized or reach consensus.
There are also adjacent problems, such as the problem where a group of spacecraft shall follow a leader while synchronizing the rotations between each other \cite{kang2002co,dimarogonas2009leader}. In the recent work by 
D. Lee \emph{et al.} \cite{lee2014asymptotic}, a dynamic level control 
scheme is presented for spacecraft formation flying with collision
avoidance.

In this work we propose six kinematic control laws.
The first two are constructed as elements of $se(3)$; 
they are linear functions of
the transformations or the relative transformations
of neighboring agents. The third and fourth are 
defined for the rotations only. They are constructed for
the tangent space of $SO(3)$ using the angular velocity. Finally, the fifth and the sixth control laws are
defined for the translations only. They are constructed in the 
tangent space of $\mathbb{R}^3$.
 All the control laws lead to 
consensus (or equivalently formation) under different assumptions on the graphs,
the initial conditions and measurable entities. 

The results for consensus on $SO(3)$ expands
on the publications \cite{distthunberg,thunberg_w,thunbergaut}, by considering 
a larger class of local representations. 
Moreover, Proposition~\ref{proposition:exponential}
provides the result that for certain topologies
and all the considered local representations, the rate of convergence is exponential. 
An interesting geometric insight is provided in 
Theorem~\ref{thm:3} where it is shown that for any 
of the local representations considered, if the
second rotation control law is used and the rotations 
initially are contained inside the injectivity region, there is a change of coordinates so that 
the system has a well known-structure.

Towards the end of this paper we also consider the second order dynamics and 
torque control laws for rigid bodies in space. We use methods similar to backstepping in order to generalize the kinematic control laws to this scenario. This generalization is only performed 
for the case of time-invariant
topologies.

The paper proceeds as follows. In Section~\ref{sec:preliminaries},
preliminary concepts are defined such as \emph{Euclidean transformations},
\emph{rotations}, \emph{translations}, \emph{network topologies} and \emph{switching signal functions}. The concept of \emph{local representations} for the rotations is also 
are introduced. In Section~\ref{chapter2:sec:problem}, the problem formulation
is given. Section~\ref{sec:kinematic:control:laws}
introduces the six kinematic control laws, which are categorized into 
two groups. Convergence results for the first group of control laws
are provided
in Section~\ref{sec:result1}, whereas the second group is treated in Section~\ref{sec:result2}. In Section~\ref{sec:dyn} -- for the 
application of rigid bodies in space -- we
provide results for control laws on the dynamic level.

\section{Preliminaries}\label{sec:preliminaries}

\subsection{Euclidean transformations, rotations, and translations}
We consider a system of $n$ agents with states in $SE(3)$, the group of Euclidean
transformations. This means that each agent $i$ has a 
matrix 
$$G_i(t) = \begin{bmatrix}
R_i(t) & T_i(t) \\
0 & 1
\end{bmatrix} \in SE(3)$$
at each time $t \geq t_0$.
The matrix $R_i(t)$ is an element of $SO(3)$, the matrix group which is 
defined by
$$SO(3) = \{R \in \mathbb{R}^{3 \times 3}: R^TR = I,  \text{det}(R) = 1\}.$$
The vector $T_i(t)$ is an element in $\mathbb{R}^3$. 

Each agent has a corresponding rigid body.
We denote the world coordinate frame by $\mathcal{F}_W$ and the instantaneous body coordinate
frame of the rigid body of each agent $i$ by $\mathcal{F}_i$. Let $R_i(t) \in SO(3)$
be the
rotation of $\mathcal{F}_i$ in the world frame $\mathcal{F}_W$ at
time $t$ and let $R_{ij}(t) \in SO(3)$ be the rotation of
$\mathcal{F}_j$ in the frame $\mathcal{F}_i$, \emph{i.e.},
$$R_{ij}(t) = R_i^T(t)R_{j}(t).$$ 
We refer to $R_i(t)$ as
\emph{absolute rotation} and $R_{ij}(t)$ as \emph{relative rotation}. 

The vector $T_i(t)$ is the position of agent $i$ in $\mathcal{F}_W$
at time $t$. The relative positions between agent $i$ and agent $j$
in the frame $\mathcal{F}_i$  at time $t$ is $$T_{ij}(t) = R_i^T(t)(T_j(t) - T_i(t)),$$
which in general is different from $T_j(t) - T_i(t)$, the relative
positions between agent $i$ and agent $j$ in the world frame. 
In the same way as for the rotations, we refer to $T_i(t)$ as 
\emph{absolute translation} and $T_{ij}(t)$ as \emph{relative translation}.

The \emph{relative Euclidean transformation} 
\begin{align*}
& G_{ij}(t) = G_i^{-1}(t)G_j(t)   \\ 
& =  \begin{bmatrix}
R_i^T(t)R_j(t) & R_i^T(t)(T_j(t) - T_i(t)) \\
0 & 1
\end{bmatrix} 
\end{align*}
contains both the relative rotation and the 
relative translation. 
From now on, in general we suppress the explicit time-dependence 
for the variables, \emph{i.e.}, $G_i$ should be 
interpreted as $G_i(t)$.

\subsection{Local representations for the rotations}\label{sec:local}
For a vector $p = [p_1, p_2, p_3]^T$ in $\mathbb{R}^3$ we define
$\widehat{{p}} = p^{\wedge}$ by
\begin{equation}\label{eq:hatmap}
\widehat{{p}} = p^{\wedge} = 
\begin{bmatrix}
  0 & -p_3 & p_2\\
  p_3 & 0 &-p_1\\
  -p_2 &p_1 & 0
\end{bmatrix}.
\end{equation}
We also define $(\cdot)^{\vee}$ as the inverse of $(\cdot)^{\wedge}$, \emph{i.e.},
 $(p^{\wedge})^{\vee} = p$.

We consider local representations or parameterizations of $SO(3)$.
Often we simply refer to them as representations or parameterizations. 
In this context, what is meant by a local 
representation is a diffeomorphism $f: B_r(I) \rightarrow B_{r', 3}(0) \subset \mathbb{R}^3$,
where $B_r(I)$ is an open geodesic ball around the identity matrix in $SO(3)$ of 
radius $r$ less than or equal to $\pi$, and $B_{r', 3}(0)$ is an open ball around the point $0$ in $\mathbb{R}^3$ with 
radius $r'$. $\bar{B}_{r}(I)$ and $\bar{B}_{r',3}(0)$ are the closures
of said balls. If we write $B_{r,3}$ or $B_r$, this is short 
hand notation for $B_{r,3}(0)$ or $B_r(I)$ respectively. The same goes for the closed balls. The local representations can be seen as coordinates 
in a chart covering an open ball around the identity matrix in $SO(3)$.

A set in $SO(3)$ is convex if any geodesic shortest path segment between any two points in the set
is contained in the set. The set is strongly convex if there is a unique geodesic shortest path segment contained in the set~\cite{afsari2011riemannian}. 
If $r = \pi$, $B_r(I)$
comprises almost all of $SO(3)$ (in terms of measure), and $B_r(I)$ is convex if and only if
$r \leq \pi/2$. The radius $r$ is referred to as the radius of injectivity. The parameterizations that we use have the following
special structure
\begin{equation}
f(R) = g(\theta)u,
\end{equation}
where $\theta$ is the geodesic distance between $I$ and 
$R$ on $SO(3)$, also referred to as the Riemannian distance,
written as $d(I,R)$. The variable $u \in \mathbb{S}^2$ is the rotational axis
of $R$, and $g: (-r, r) \rightarrow \mathbb{R}$ is an odd, analytic and
strictly increasing function such that $f$ is a diffeomorphism.
On $B_{\pi}(I)$ the vector $u$ and the positive variable $\theta$ are obtained
as functions of $R$ in the following way
$$\theta = \cos^{-1}\bigg ( \frac{\text{trace}(R) -1}{2}\bigg ), \quad u = \frac{1}{2\sin(\theta)}
\begin{bmatrix}
r_{32} - r_{23} \\
r_{13} - r_{31} \\
r_{21} - r_{12}
\end{bmatrix},$$
where $R = [r_{ij}]$. 
Let us denote $y_i = f(R_i)$  and $y_{ij} = f(R_{ij}).$
It holds that
$
{y}_{ij}  =  -{y}_{ji}, \text{ but in general } {y}_{j} - {y}_{i} \neq {y}_{ij}.$
For each representation, \emph{i.e.}, choice of $g$, $r \leq \pi$
is the largest radius such that $f$ is a diffeomorphism.
The radius $r$ is the radius of injectivity and depends on the representation, but we suppress this explicit dependence and
throughout this paper, 
$r$ corresponds to the representation at hand, \emph{i.e.}, the one we have chosen to consider at the moment. For the representation at hand we also define 
$$r' = \sup_{s \uparrow r}g(s).$$
Some common representations are:
\begin{itemize}
\item \textbf{ The Axis-Angle Representation}, in which case $g(\theta) = \theta$ and
$r = r' =\pi$. This representation is almost global. The set $SO(3)\backslash B_{\pi}(I)$ has measure zero in $SO(3)$. 
The Axis-Angle Representation is obtained from the
logarithmic map by
\begin{align*}
{{x}}_i& = (\text{Log}(R_i))^{\vee}, \\
{{x}}_{ij}& = (\text{Log}(R_i^TR_j))^{\vee}.
\end{align*}
In the other direction, a rotation matrix $R_i$ is obtained via the exponential 
map by 
$$R_i(x_i) = \text{exp}(\widehat{x}_i).$$
The matrix $R_{ij}$  is obtained by 
$$R_{ij}(x_i, x_j) =  \text{exp}(\widehat{x}_i)^T\text{exp}(\widehat{x}_j).$$
The function $ \text{exp}_{R_i}$ is the exponential map at $R_i$. Using this notation, the function $\text{exp}$ is short hand notation for $\text{exp}_{I}$.
\\

\item \textbf{ The Rodrigues Parameters}, in which case $g(\theta) = \tan(\theta/2)$. The corresponding 
$r$ and $r'$ are equal to $\pi$ and $\infty$ respectively. \\
\item \textbf{ The Modified Rodrigues Parameters}, in which case $g(\theta) = \tan(\theta/4)$, $r = \pi$ and $r' = 1$.
This representation is obtained from the rotation matrices by a second order Cayley transform~\cite{tsiotras1997higher}.\\
\item \textbf{ The representation $(R -R^T)^{\vee}$}, in which case $g(\theta) = \sin(\theta)$, 
and the corresponding $r$ and $r'$ are $\pi/2$ and $1$ respectively. This representation is popular because it is easy to express
in terms of the rotation matrices. Unfortunately, since $r = \pi/2$, only $B_{\pi/2}(I)$ is covered. \\
\item \textbf{ The Unit Quaternions}, or rather parts of it. The unit quaternion $q_i$, expressed as a function
of the Axis-Angle Representation $x_i = \theta_iu_i$ of $R_i \in B_{\pi}(I)$, is given by
$$q(x_i) = (\cos(\theta_i/2), \sin(\theta_i/2)u_i)^T \in \mathbb{S}^3.$$
This means that we can choose the last three elements of the unit
quaternion vector as our representation, \emph{i.e.}, $\sin(\theta_i/2)u_i$,
in which case $r = \pi$. 
The unit quaternion representation is popular since the mapping from $SO(3)$ 
to the quaternion sphere is a Lie group homomorphism.
\end{itemize}

Let ${x}_i(t)$ and ${x}_{ij}(t)$ denote the axis-angle
representations of the rotations $R_i(t)$ and $R_{ij}(t)$,
respectively. In the following, since we are only addressing 
representations of (subsets of) $B_{\pi}(I)$, we choose ${x}(t) = [{x}_1^T(t), {x}_2^T(t), \ldots, {x}_n^T(t)]^T \in (B_{\pi,3}(0))^n$
as the state of the system instead of $(R_1(t), \ldots, R_n(t)) \in (B_{\pi}(I))^n.$ 
Note that since $\theta_i = \|x_i\|$, it holds that $g(\theta_i) = g(\|x_i\|)$. The variables $y_i$ and $y_{ij}$ can be seen as functions of $x_i$ and $x_i, x_j$
respectively,
\emph{i.e.}, 
\begin{align*}
y_i(x_i) & = (f \circ \text{exp})(\widehat{x}_{i}), \\
y_{ij}(x_i, x_j) & = (f \circ \text{exp})(\text{Log}(R_i^T(x_i)R_j(x_j))).
\end{align*}

Since $x_i$ and $x_j$ are elements of the vector $x$, we can write $y_i(x)$ and $y_{ij}(x)$.
When we write $y_i(t)$ and $y_{ij}(t)$, this 
is equivalent to $y_{i}(x(t))$ and $y_{ij}(x(t))$ respectively. 
If we want to emphasize the dependence of the initial condition,
instead of writing $x(t)$ (or $y(t)$) we write $x(t,t_0,x_0)$
(or $y(t,t_0,y_0))$ where $x_0$ is the initial state and 
$t_0$ is the initial time.

\subsection{Kinematics}
We denote the instantaneous angular velocity of $\mathcal{F}_i$
by ${\omega}_i$. From now on, until Section~\ref{sec:dyn}, we assume that $\omega_i$ is the 
control variable for the rotation of agent $i$. 
The kinematics for $R_i$ is given by 
$$\dot{R}_i = R_i\widehat{\omega}_i,$$
where $R_i\widehat{\omega}_i$ is an element of the tangent space $T_{R_i}SO(3)$.

The kinematics
is given by
\begin{align}\label{dynamics}
\dot{{x}}_{i} & = L_{{x}_{i}}{\omega}_{i},
\end{align}
 where the Jacobian (or transition) matrix $L_{{x}_i}$ is given by
\begin{equation}\label{jacobian}
L_{{x}_i} = L_{\theta_i {u}_i} = I_3 +
\frac{\theta}{2}\widehat{{u}}_i+ \bigg (1 -
\frac{\text{sinc}(\theta_i)}{\text{sinc}^2(\frac{\theta_i}{2})}\bigg
)\widehat{{u}}_i^2.
\end{equation}
The proof is found in \cite{junkins2003analytical}. 
The function $\mathrm{sinc}(\beta)$ is
defined so that $\beta\mspace{3mu}\mathrm{sinc}(\beta) = \sin(\beta)$ and $\mathrm{sinc}(0) =
1$. It was shown in \cite{malis19992} that
$L_{\theta {u}}$ is invertible for $\theta \in
(-2\pi,2\pi)$. Note however that $\theta \in [0, \pi)$ here.

The linear velocity of agent $i$, expressed in $\mathcal{F}_i$, is 
denoted by $v_i$. Up until Section~\ref{sec:dyn},
we assume that $v_i$ is the control variable 
for the translation of agent $i$. The time derivative of $T_i(t)$ is given by 
$$\dot{T}_i(t) = R_i(t)v_i.$$
Define
$$\xi_i = \begin{bmatrix}
\widehat{\omega}_i & v_i \\
0 & 0
\end{bmatrix}.$$
It holds that 
$$\dot{G}_i(t) = G_i(t)
\xi_i.$$
\subsection{Dynamics}
The dynamics for agent $i$ is given by
\begin{align*}
\dot{G}_i & = G_i\xi_i \\
\dot{\xi}_i & = 
\begin{bmatrix}
(J_i^{-1}(-\widehat{\omega}_iJ_i\omega_i + \boldsymbol{\tau}_i))^{\wedge} & (\frac{\boldsymbol{f}_i}{m_i} - \widehat{\omega}_iv_i) \\
0 & 0
\end{bmatrix},
\end{align*}
where $J_i$ is the inertia matrix, $m_i$ is the mass, $\boldsymbol{\tau_i}$ is the control torque, and $\boldsymbol{f}_i$ is control force -- the latter two are given as a bold symbols
since we do not want to mix them up with other defined entities.

\subsection{Connectivity}
\begin{defn}\label{def:chapter1:graff}
A directed graph (or digraph) $\mathcal{G} = (\mathcal{V},
\mathcal{E})$ consists of a set of nodes,  $\mathcal{V} = \{1, ..., n\}$ 
and a set of edges $\mathcal{E} \subset \mathcal{V} \times \mathcal{V}$.
\end{defn}
Each node in the graph corresponds to a unique agent.
We also define neighbor sets or neighborhoods.
Let $\mathcal{N}_i \in
\mathcal{V}$ comprise the neighbor set (sometimes referred to simply as neighbors)
of agent $i$, where $j \in \mathcal{N}_i$
if and only if $(i,j) \in \mathcal{E}$. We assume that $i \in
\mathcal{N}_i$ \emph{i.e.}, we restrict the collection of graphs
to those for which $(i,i) \in \mathcal{E}$ for all $i \in \mathcal{V}$.

A directed path of $\mathcal{G}$ is an ordered sequence of distinct nodes in $\mathcal{V}$
such that any consecutive pair of nodes in the sequence corresponds to an edge in the graph.
An agent $i$ is connected to an agent $j$ if there is a directed path starting in $i$ and ending in $j$.

\begin{defn}
A digraph is strongly connected if each node $i$
is connected to all other nodes.
\end{defn} 

\begin{defn}
A digraph is quasi-strongly connected if there exists a rooted 
spanning tree or a center, {i.e.}, at least one node such
that all other nodes are connected to it.
\end{defn}

An adjacency matrix $\mathcal{A} = [a_{ij}]$ for a graph
$\mathcal{G} = (\mathcal{V},
\mathcal{E})$ is a matrix where $a_{ij} \geq 0$ for all $i,j$, and furthermore $a_{ij} > 0$
if and only if $(i,j) \in \mathcal{E}$ for all $i,j$. Given our definition of graph, \emph{i.e.}, Definition~\ref{def:chapter1:graff},
there are infinitely many adjacency matrices for a graph.

From Definition~\ref{def:chapter1:graff} we see that there are $2^{n^2}$
directed graphs with $n$ nodes, \emph{i.e.}, the power set of the 
edge set. Since we assume that 
$(i,i)$ is an edge in the graph for all $i$, there are $2^{n^2 - n}$ graphs we consider. 
For $k \in \{1, \ldots, 2^{n^2 - n}\}$ we associate a corresponding unique graph $\mathcal{G}_k = (\mathcal{V}, \mathcal{E}_{k})$ and a 
unique adjacency matrix ${A}_k$.
The ${A}_k$ matrices are constructed in the following way.
We construct a positive adjacency matrix ${A}' = [a_{ij}]$ for the complete (fully connected) graph. For the matrix ${A}_k = [a_{ij}^k]$,
it holds that $a_{ij}^k = a_{ij}$ if $(i,j) \in \mathcal{E}_k$,
otherwise $a_{ij}^k = 0$. Thus, if $(i,j) \in \mathcal{E}_k$,
we can write $a_{ij}$ instead of $a_{ij}^k$.

Now, for each agent $i$ there are $2^{n-1}$ unique neighborhoods 
$\mathcal{N}_i^l$, where $l \in \{1, \ldots, 2^{n-1}\}$.  
Given $k \in \{1, \ldots, 2^{n^2 - n}\}$, for agent $i$
there is a unique $l \in \{1, \ldots, 2^{n - 1}\}$ such that 
$\mathcal{N}_i^l$ is the neighborhood of agent $i$ in the graph $\mathcal{G}_k$.
Also, if each agent $i$ has chosen an $l \in \{1, \ldots, 2^{n - 1}\}$
such that $\mathcal{N}_i^l$ is the neighborhood of agent $i$,
then there is a unique $k \in \{1, \ldots, 2^{n^2 - n}\}$ such that
$\mathcal{G}_k$ is the graph for the system.

We are now ready to address time-varying graphs.
In order to do so, for each agent $i$, we introduce 
a switching signal function 
$$\sigma^i: \mathbb{R} \rightarrow \{1, \ldots, 2^{n-1}\},$$
which is piece-wise constant and right-continuous.
Let $\{\tau^i_k\}$ be the monotonically strictly increasing  
sequence of times for which $\sigma^i$
is discontinuous. 
We assume that there is a positive lower bound $\tau_D$ between two consecutive switches, \emph{i.e.}, 
$$\sup_{k}(\tau^i_{k+1} - \tau^i_k) > \tau_D \quad \text{ for all } i.$$
The time-varying neighborhood of agent $i$ is
$\mathcal{N}^{\sigma^i(t)}_i$. 

Given the set of switching signal functions $\sigma^i$ we 
can construct a piece-wise constant and right-continuous 
switching signal function for the graph of the
multi-agent system. This switching signal function $\sigma$
has range $\{1, \ldots, 2^{n^2 - n}\}$ and 
switching times
$$\{\tau_k\} = \bigcup_{i}\{\tau_{l}^i\},$$
where $\{\tau_k\}$ is monotonically strictly increasing in $k$. Note that
for $\sigma$ it is not 
necessarily true that there is a positive lower bound on the dwell time
between two consecutive switches as is the case for $\sigma^i$. 

Now,
between any two switching times,
$\sigma(t)$ is equal to the $k$ for which the graph $\mathcal{G}_k$ it holds that 
the neighborhood of each agent $i$ is equal to $\mathcal{N}_i^{\sigma^i(t)}$.

\begin{defn}
The union graph of $\mathcal{G}_{\sigma(t)}$ during the time
interval $[t_1,t_2)$ is defined by
\begin{equation*}
\mathcal{G}([t_1, t_2))
= \textstyle\bigcup_{t\in[t_1, t_2)} \mathcal{G}_{\sigma(t)}
= (\mathcal{V},\textstyle\bigcup\nolimits_{t\in[t_1, t_2)}\mathcal{E}_{\sigma(t)}),
\end{equation*}
where $t_1 < t_2 \leq +\infty$.
\end{defn}

\begin{defn}\label{def:quasi}
The graph $\mathcal{G}_{\sigma(t)}$ is uniformly
(quasi-) strongly connected if there is $T^{\sigma}>0$ such that the
union graph $\mathcal{G}([t, t + T^{\sigma}))$ is (quasi-) strongly connected for all $t$.
\end{defn}

The idea of using an individual switching signal $\sigma^i$ for each agent, is that each agent shall be able to choose independently 
which neighbors it decides to receive information from.

Instead of using the term \emph{communication graph} for
$\mathcal{G}_{\sigma(t)}$, we deliberately use the terms
\emph{neighborhood graph}, \emph{connectivity graph} or \emph{interaction graph}. Direct communication does not necessarily take place
between the agents in practice. Instead, they can choose to just observe each
other via cameras or other sensors, \emph{i.e.}, indirect communication.

\section{Consensus and formation control}\label{chapter2:sec:problem}
\subsection{Consensus}
We start this section by introducing the consensus
problem on $SE(3)$. Consensus on $SE(3)$
means that, as time tends to infinity, 
the set of transformations $(G_1(t), G_2(t), \ldots, G_n(t)) \in (SE(3))^n$ approaches 
the consensus set where all the transformations are
equal. The problem is to construct a distributed control law 
for each agent $i$, where only information from the neighbors $\mathcal{N}_i$ is used in the control law, such that the system reaches consensus. This information could be the relative transformations to the neighbors or the absolute transformations of the neighbors. 
An other desired property is that 
the velocities $\xi_i$ tend to zero sufficiently fast so that the
transformations converge to a static transformation. 

When we say that the Euclidean transformations of the 
agents \enquote{approaches} the consensus set, we mean that 
the rotations $$(R_1(t), R_2(t), \ldots, R_n(t)) \in SO(3)^n$$
approach $\{(R_1, \ldots, R_n) \in (\bar{B}_{q}(I))^n: R_1 = \hdots = R_n\}$ and the translations 
$$T^{\text{tot}} = [T_1^T(t), T^T_2(t), \ldots, T^T_n(t)]^T \in \mathbb{R}^{3n}$$
 approach the set
where all the translations are equal. For the translations 
the convergence is defined in terms of the Euclidean metric.
For the rotations, the convergence is defined in terms of the
Riemannian metric on $SO(3)$. If the rotations are contained within the region of injectivity of a local parameterization, asymptotic stability in terms of the Riemannian metric on $SO(3)$ is equivalent to asymptotic stability using the Euclidean metric in the parameterization domain for $x$.

The 
consensus problem on $SE(3)$ might seem uninteresting in practice,
since for rigid bodies in space it is not physically possible 
to reach consensus in the positions. There are two reasons for
considering this problem anyway. Firstly, if we look at the consensus
problem as two subproblems, consensus in the rotations and consensus
in the positions, the former is still interesting in practice and
has received a great deal of attention lately. Secondly and more importantly, the consensus control problem is equivalent to the formation 
control problem after a change of coordinates. Thus, all the control laws we develop for the consensus
control problems can also be used for the formation control problem after a simple
transformation. This will be elaborated more in Section~\ref{sec:form}.

The subproblem of reaching consensus in the rotations is referred to as the \emph{attitude synchronization} problem or \emph{consensus on $SO(3)$}. Then we shall find a feedback control law ${\omega}_i$ for each agent $i$ using the local representations of either absolute rotations or relative rotations so that the absolute rotations of all agents converge to the set where all the rotations are equal as time goes to infinity, \emph{i.e.},
\begin{equation}\label{condition1}
\|R_i(t) - R_j(t)\| \rightarrow 0,    \text{ for all } i,j, \text{ as } t \rightarrow \infty,
\end{equation}
or equivalently,
\begin{equation*}
\hspace{9mm} \|R_{ij}(t) - I\| \rightarrow 0,  \text{ for all } i,j,    \:\: \text{as } t \rightarrow \infty.
\end{equation*}
If $y \in (B_{r',3}(0))^n$ it is true that 
\begin{align}\label{eq:z1}
R_i =R_j & \Longleftrightarrow x_i = x_j \Longleftrightarrow x_{ij} = 0  \\
\nonumber
& \Longleftrightarrow y_i = y_j  \hspace{0.51mm} \Longleftrightarrow y_{ij} = 0 \quad \text{ for all } i,j.  
\end{align}
We define the consensus set $\mathcal{A}$ in $\mathbb{R}^{3n}$ as follows:
$$\mathcal{A} = \{z = [z_1^T, z_2^T, \ldots, z_n^T]^T \in \mathbb{R}^{3n}: z_i = z_j \in \mathbb{R}^3, \forall i,j\}.$$
According to \eqref{eq:z1} and the fact that the 
map $$R_i \mapsto x_i$$ is a diffeomorphism on $B_{\pi}(I)$,  \eqref{condition1} can equivalently be 
written as $x(t) \rightarrow \mathcal{A}$ as $t\rightarrow \infty$.
This means that the solution approaches $\mathcal{A}$.
Thus, provided we can guarantee that $y(t) \in (B_{r',3}(0))^n$
for all $t \geq t_0$, where $t_0$ is the initial time,
consensus on $SE(3)$ for the multi-agent system is the following
$$(x(t),T^{\text{tot}}(t)) \rightarrow \mathcal{A} \times \mathcal{A}, \quad \text{ as } t \rightarrow \infty. $$

A stronger assumption on the convergence to $\mathcal{A} \times \mathcal{A}$ is global uniform asymptotic stability of  $\mathcal{A} \times \mathcal{A}$
relative to a strongly forward-invariant set, see Definition \ref{def:strongly:forward:invariant} and Definition~\ref{def:uniform:stability} below. The distance from 
a point $z$ in $\mathbb{R}^p$ to a set $\mathcal{D}$ in $\mathbb{R}^p$ is defined by
$$\|z\|_{\mathcal{D}} = \inf_{w \in \mathcal{D}}\|z-w\|.$$

For a time-invariant system, forward invariance or positive invariance  of a set means that every solution to the system with initial condition
in the set is forward complete and the solution at any time greater than the initial time is contained in the set. For switched systems we have the following type of invariance.
The $f$-vectors used in the following two definitions are locally 
defined in that context.

\begin{defn}\label{def:strongly:forward:invariant}
Consider dynamical systems of the following class. The 
dynamical equation is given by 
$$\dot{z} = f_{\sigma(t)}(z),$$
where $z(t) \in \mathbb{R}^p$ for some positive integer $p$.
The right-hand side is switching between a finite set $\mathcal{F} = \{f_k\}$
of time-invariant functions according 
to a switching signal function $\sigma$. The switching signal function $\sigma$ is well-behaved in the sense that there are only 
finitely many switches on any compact time interval. 

A set $\mathcal{D} \subset \mathbb{R}^p$ is strongly forward-invariant if 
for any time $t_0$, any $z_0 \in \mathcal{D}$ and any such well behaved switching signal function $\sigma$ switching between functions in $\mathcal{F}$, the solution $z(t,t_0,z_0)$
exists, is unique, forward complete and contained in $\mathcal{D}$ for all $t \geq t_0$.
\end{defn}

\begin{defn}\label{def:uniform:stability}
Consider the dynamical system 
$$\dot{z} = f(t,z),$$
where $y(t) \in \mathbb{R}^p$ for some positive integer $p$.
A set $\mathcal{D}_1 \subset \mathbb{R}^p$ is globally
uniformly asymptotically stable relative to the compact
strongly forward-invariant set $\mathcal{D}_2$, if 
\begin{enumerate}
\item $\mathcal{D}_1$ is uniformly stable relative to $\mathcal{D}_2$,
{i.e.}, for every $\epsilon > 0$, there is a $\delta(\epsilon) > 0$ such that 
\begin{align*}
& (\|z_0\|_{\mathcal{D}_1} \leq \delta, \: z_0 \in \mathcal{D}_2) \Longrightarrow \\
& (\|z(t_2,t_1,z_0)\|_{\mathcal{D}_1} \leq \epsilon \: \textnormal{ for all } t_1, t_2 \textnormal{ where } t_2 \geq t_1),
\end{align*}
\item $\mathcal{D}_1$ is globally uniformly attractive relative to 
$\mathcal{D}_2$, {i.e.},
for every $\epsilon > 0$, there is a $\tau(\epsilon) > 0$ such that
\begin{align*}
& z_0 \in \mathcal{D}_2 \Longrightarrow \\
& (\|z(t_2,t_1,z_0)\|_{\mathcal{D}_1} \leq \epsilon \: \textnormal{ for all } t_1,t_2, \\
& \textnormal{ such that } t_2 \geq t_1 + \tau(\epsilon)).
\end{align*} 
\end{enumerate}
\end{defn}

One can show that if $\mathcal{A}$ is globally uniformly
asymptotically stable relative to the strongly forward invariant set $(\bar{B}_{q,3}(0))^n$ for $x(t)$ where $q < \pi$, then the set $\{(R_1, \ldots, R_n) \in (\bar{B}_{q}(I))^n: R_1 = \hdots = R_n\}$ is globally
uniformly asymptotically stable relative to $(\bar{B}_{q}(I))^n$  (when the Riemannian metric is used). The notation of \enquote{strong forward invariance} is adopted from~\cite{goebel2012hybrid},
where it is defined for hybrid systems.

\subsection{Formation}\label{sec:form}
The consensus problem has many applications in the cases where the motion
is purely rotational, \emph{e.g.}, attitude synchronization for spacecraft or
orientation alignment for cameras. However, as already mentioned,  reaching consensus in the 
positions is obviously not physically possible for rigid bodies, but reaching 
a formation is.

The objective is to make the $G_i^{-1}(t)G_j(t)$ matrices 
converge to some desired $G^*_{ij}$ matrices. The $G^*_{ij}$
matrices are assumed to be transitively consistent in that 
$$G^*_{ij}G^*_{jk} = G^*_{ik} \text{ for all } i,j,k.$$
A necessary and sufficient condition for transitive consistency~\cite{tron2014distributed,bernard2014,thunberg2015transitive}
of the $G^*_{ij}$ is that there are $G^*_{i}$ such that 
$$G^*_{ij} = G^{*-1}_{i}G^*_{j} \text{ for all }i,j.$$
In this light, we formulate the objective in the formation problem as
follows. Given some desired constant Euclidean transformation matrices $G_1^*, \ldots, G_n^*$,
construct a control law for each agent $i$ such that
\begin{align*}
~\|G_1\left(t\right) - Q^{-1}(t)G_1^*\| & \rightarrow 0, \\
~\|G_2\left(t\right) - Q^{-1}(t)G_2^*\| & \rightarrow 0, \\
~& \: \vdots \\
~\|G_n\left(t\right) - Q^{-1}(t)G_n^*\| & \rightarrow 0, 
\end{align*}
as $t \rightarrow \infty$, where $Q(t)$ is a Euclidean transformation. This implies that
\begin{align*}
\|G_i^{-1}\left(t\right)G_j\left(t\right) - G_i^{*-1}G_j^*\| \rightarrow 0 \quad \text{ as } \quad  t \rightarrow \infty.
\end{align*}
Thus, in some (possibly time-varying) coordinate frame,
the Euclidean transformation of agent $i$ converges to $G_i^*$
as time tends to infinity.
Each matrix $G_i^*$ contains the rotation matrix $R_i^*$ and the 
translation $T_i^*$.

On a kinematic level the formation control problem is 
equivalent to the consensus problem. Let us define 
$$\tilde{G}_i = G_iG_i^{*-1} \: \text{ and } \: \tilde{G}_{ij} = G_i^{*}G_{ij}G_j^{*-1}, \quad \text{for all } i.$$
The kinematics for $\tilde{G}_i$ is given by
$$\dot{\tilde{G}}_i = G_i\xi_iG_i^{*-1} = \tilde{G}_iG_i^{*}\xi_iG_i^{*-1}
= \tilde{G}_i\tilde{\xi}_i,$$
where 
\begin{align*}
 \tilde{\xi}_i &  = G_i^{*}\xi_iG_i^{*-1} = \begin{bmatrix}
\widehat{\tilde{\omega}}_i & \tilde{v}_i\\
0 & 0
\end{bmatrix}  \\
& = 
\begin{bmatrix}
R_i^*\widehat{\omega}_iR_i^{*T} & - R_i^*\widehat{\omega}_iR_i^{*T}T_i^* + R_i^*v_i\\
0 & 0
\end{bmatrix}.
\end{align*}
and 
\begin{align*}
 {\xi}_i &  = G_i^{*-1}\tilde{\xi}_iG_i^* = \begin{bmatrix}
\widehat{{\omega}}_i & {v}_i\\
0 & 0
\end{bmatrix}  \\
& = 
\begin{bmatrix}
R_i^{*T}\widehat{\tilde{\omega}}_iR_i^{*} &  R_i^{*T}(\widehat{\tilde{\omega}}_iT_i^* + \tilde{v}_i)\\
0 & 0
\end{bmatrix}.
\end{align*}
It easy to see that if the system reaches consensus in the $\tilde{{G}}_i$, it also reaches the desired formation.
Thus, a consensus control law $\tilde{\xi}_i$ can be constructed
for each agent and provided that each agent $i$ knows $G_i^*$,
$\xi_i$ is obtained by $\xi_i = G_i^{*-1}\tilde{\xi}_iG_i^{*}$.
In general, unless the design is 
limited to the $\omega_i$, the proposed control laws in the next section
should be used for formation control of the $\tilde{\xi}_i$.

On the dynamic level we have that
\begin{align}
\label{eq:formation_dyn_1}
\dot{\tilde{\omega}}_i  =~ & R_i^{*}J_i^{-1}\left(-\left(R_i^{*T}\widehat{\tilde{\omega}}_i R_i^{*}\right)J_iR_i^{*T}\tilde{\omega}_i + \boldsymbol{\tau}_i\right), \\
\label{eq:formation_dyn_2}
\dot{\tilde{v}}_i  =~ & R_i^{*}\left (\frac{\boldsymbol{f}_i}{m_i} - R_i^{*T} \widehat{\tilde{\omega}}_i\left (\widehat{\tilde{\omega}}_iT_i^* + \tilde{v}_i\right ) \right ) - \left (\dot{\tilde{\omega}}_i\right)^{\wedge} T_i^*.
\end{align} 
For the control design on the dynamic level, the approach is to 
design a consensus control law for the $\tilde{\omega}_i$ and the
$\tilde{v}_i$ and then track 
this desired kinematic control law using methods similar to backstepping. 
For control laws designed on the kinematic level, since the 
problems of consensus and formation are equivalent, we will only 
focus on the consensus problem. The consensus problem more 
tractable, since one can use existing theory for that problem. On the dynamic level we will also only 
consider the consensus problem -- the formation control
laws have a similar structure as the consensus control laws 
in this case.

\section{Kinematic control laws}\label{sec:kinematic:control:laws}
We use two approaches for the design of the $\xi_i$. 
The first approach is to treat $\xi_i$ as one control variable 
and design a feedback control law as an expression of the $G_i$,
the second approach is to design $\omega_i$ and $v_i$ separately. 
Most emphasis will be on the second approach. The control laws in the first 
approach are referred to as the \emph{first control laws},
whereas the control laws in the second approach are referred to as \emph{the second control laws}.

\subsection*{The first control laws}
We propose the following control laws based on absolute and relative transformations respectively.
\begin{align}
\label{eq:approach1:controller:1}
{\xi}_i  & = \sum_{j\in\mathcal{N}^{\sigma^i(t)}_i} a_{ij}\left( \left(G_j - G_i) + (G_i^{-1} - G_j^{-1}\right) \right), \\
\label{eq:approach1:controller:2}
{\xi}_i  & = \sum_{j\in\mathcal{N}^{\sigma^i(t)}_i} a_{ij}\left( G_{ij} - G_{ij}^{-1} \right).
\end{align}

\subsubsection*{The second control laws}
In the first two control laws below,
$y_i$ and $y_{ij}$ could be any of the local representations considered in Section~\ref{sec:preliminaries}.
\begin{align}
\label{chapter2:controller:1}
{\omega}_i  & = \sum_{j\in\mathcal{N}^{\sigma^i(t)}_i} a_{ij}({y}_{j} - {y}_{i}), \hspace{29mm}\\
\label{chapter2:controller:2}
{\omega}_i &= \sum_{j\in\mathcal{N}^{\sigma^i(t)}_i} a_{ij}{y}_{ij}, \\
\label{chapter2:controller:1:v}
{v}_i  & = \sum_{j\in\mathcal{N}^{\sigma^i(t)}_i} a_{ij}({T}_{j} - {T}_{i}), \hspace{29mm}\\
\label{chapter2:controller:2:v}
{v}_i &= \sum_{j\in\mathcal{N}^{\sigma^i(t)}_i} a_{ij}{T}_{ij}.
\end{align}
The structure of these second control laws and especially 
\eqref{chapter2:controller:1} and \eqref{chapter2:controller:1:v}
are well known from the 
literature \cite{mesbahi2010graph,ren2007distributed}. In Section \ref{sec:result2} we provide
new results on the rate of convergence and regions of attractions for these 
control laws in this context. When the control laws are used for 
formation instead of consensus, the $\tilde{\xi}_i$ are designed instead of 
the $\xi_i$; the controllers are obtained through the relation 
$$\tilde{\xi}_i   = G_i^{*}\xi_iG_i^{*-1},$$ 
as given in Section~\ref{sec:form}. As an example, suppose all the $R_i$ rotations and all 
the desired $R_i^*$ rotations in the formation are equal to the identity matrix. Then the
agents shall reach a desired formation in the positions only. All the agents 
construct $\tilde{v}_i$ according to \eqref{chapter2:controller:1:v} or \eqref{chapter2:controller:2:v} and solve for $v_i$ through the following 
relation
$$\begin{bmatrix}
I & T_i^* \\
0 & 1
\end{bmatrix}^{-1}
\begin{bmatrix}
0 & \tilde{v}_i \\
0 & 0
\end{bmatrix}
\begin{bmatrix}
I & T_i^* \\
0 & 1
\end{bmatrix}
=
\begin{bmatrix}
0 & \tilde{v}_i \\
0 & 0
\end{bmatrix}
=
\begin{bmatrix}
0 & v_i \\
0 & 0
\end{bmatrix}.
$$
In this simple case $v_ i = \tilde{v}_i$ and $\xi_i = \tilde{\xi}_i$. However, in general $\xi_i \neq \tilde{\xi}_i$.

The following two sections are devoted to the study of the control laws (\ref{eq:approach1:controller:1}-\ref{chapter2:controller:2:v}).

\section{Results for the first control laws}\label{sec:result1}

\begin{prop}\label{proposition:first_control_law:1}
Suppose the graph $\mathcal{G}_{\sigma(t)}$ is time-invariant and strongly connected.
Suppose that each rotation is contained in $B_{\pi/2}(I)$, then if control law \eqref{eq:approach1:controller:1} is used, the set $(B_{\pi/2,3}(0))^n$ is strongly forward invariant for the dynamics of $x$ and 
$$(x(t), T^{\text{tot}}(t)) \rightarrow \mathcal{A} \times \mathcal{A} \quad \text{ as } t \rightarrow \infty.$$
\end{prop}

\begin{prop}\label{proposition:first_control_law:2}
Suppose the graph $\mathcal{G}_{\sigma(t)}$ is time-invariant and
quasi-strongly connected. Suppose $q < \pi/4$, if all the rotations are contained in $\bar{B}_{q}(I)$, then if control law \eqref{eq:approach1:controller:2} is used, the set $(\bar{B}_{q,3}(0))^n$ is strongly forward invariant for the dynamics of $x$ and 
$\mathcal{A} \times \mathcal{A}$
is globally asymptotically stable relative to 
$(\bar{B}_{q,3}(0))^n \times \mathbb{R}^{3n}$.
\end{prop}

\begin{rem}
It can be shown that the results in propositions \ref{proposition:first_control_law:1} and \ref{proposition:first_control_law:2} are slightly more general. It is true 
that $(x(t), T^{\text{tot}}(t))$ converges to a fixed point in $\mathcal{A} \times \mathcal{A}$, {i.e.}, not a limit cycle. This result is however not shown here.
\end{rem}

In proposition~\ref{proposition:first_control_law:2}, we only guarantee \emph{stability} of a set instead of \emph{uniform stability} of the set. 

In the following two proofs, since the graph is time-invariant, we write 
$\mathcal{G}$ and $\mathcal{N}_i$ instead of
$\mathcal{G}_{\sigma(t)}$ and $\mathcal{N}^{\sigma^i(t)}_i$ respectively.
The graph Laplacian
matrix $L(\mathcal{G},A)$ for the graph $\mathcal{G}$ with the adjacency matrix $A$, is $$L(\mathcal{G},A) = D - A,$$
where 
 $$D = \text{diag}(d_1, \ldots, d_n) =  \text{diag}\left(\sum_{j=1}^n a_{1j}, \ldots, \sum_{j=1}^n a_{nj}\right).$$
 
\quad \emph{Proof of Proposition~\ref{proposition:first_control_law:1}:}
When the control law \eqref{eq:approach1:controller:1}
is used, $\omega_i$ is given by the following expression
$${\omega}_i  = \sum_{j\in\mathcal{N}_i} a_{ij}({y}_{j} - {y}_{i}),$$
where $y_i = \sin(\theta_i)u_i$ for all $i$. This control law for $\omega_i$ is on the form \eqref{chapter2:controller:1} and we
will later show that, provided the rotations are contained within
the region of injectivity, which in this case is the ball around
the identity with radius $\pi/2$, $x(t)$ approaches $\mathcal{A}$ asymptotically. Also,
$(B_{\pi/2,3}(0))^n$ is forward invariant.

Given the initial states $x_i(t_0)$, since there are finitely many agents, there is a positive $q < \pi/2$ such that $x(t_0) \in (\bar{B}_{q,3}(0))^n$.
Let
$\mathcal{X} = (\bar{B}_{q,3}(0))^n \times \mathbb{R}^{3n}$ and define the
two closed sets
\begin{align*}
\Gamma_2 & = \mathcal{A} \cap (\bar{B}_{q,3}(0))^n \times \mathbb{R}^{3n} \\
\Gamma_1 & = \mathcal{A} \cap (\bar{B}_{q,3}(0))^n \times \mathcal{A}.
\end{align*} 
We can choose the state space as $\mathcal{X}$ for 
$(x, T^{\text{tot}})$ since this set is 
forward invariant, see Proposition~\ref{thm:1} in Section~\ref{sec:result2}.
We observe that $\Gamma_1 \subset \Gamma_2 \subset \mathcal{X}.$

On $\mathcal{X}$, the dynamics for $T_i$ is given by 
$$\dot{T}_i = \sum_{j\in\mathcal{N}_i}a_{ij}\left(\left(R_i + R_iR_j^T)\right)T_j - \left(R_i + I\right)T_i\right).$$
But on the set $\Gamma_2$ the dynamics for $T_i$ is given by 
$$\dot{T}_i = \sum_{j\in\mathcal{N}_i}a_{ij}\left(\left(I + Q^*\right)\left(T_j - T_i\right)\right), \hspace{20mm}$$
where $Q^* \in \bar{B}_{q,3}(0)$ is some constant rotation matrix. On $\Gamma_2$,
the dynamics for $T^{\text{tot}}$ is given by 
$$\dot{T}^{\text{tot}} = -(L(\mathcal{G},A) \otimes (I + Q^*)) T^{\text{tot}}.$$

By using the fact that the eigenvalues of $(I + Q^*)$ have
real parts strictly greater than zero, the fact that ($(B_{\pi/2,3}(0))^n$ is forward invariant), 
and the fact that $L(\mathcal{G},A)$ is the graph Laplacian matrix for 
a strongly connected graph, one can show that 
$\Gamma_1$ is exponentially stable relative to $\Gamma_2$.  
Now one can use Theorem 8 in \cite{reduction} in order to show that
$\Gamma_1$ is globally attractive relative to $\mathcal{X}.$
\hfill $\blacksquare$ \vspace{3mm}

\quad \emph{Proof of Proposition~\ref{proposition:first_control_law:2}:}
When the control law \eqref{eq:approach1:controller:2}
is used, $\omega_i$ is given by the following expression
$${\omega}_i  = \sum_{j\in\mathcal{N}_i} a_{ij}y_{ij},$$
where $y_{ij} = \sin(\theta_{ij})u_{ij}$ for all $i$. This control law for $\omega_i$ is on the form \eqref{chapter2:controller:2}.

Let 
$\mathcal{X} = (\bar{B}_{q, 3}(0))^n \times \mathbb{R}^{3n}$ and define the
two closed sets
\begin{align*}
\Gamma_2 & = \mathcal{A} \cap (\bar{B}_{q, 3}(0))^n \times \mathbb{R}^{3n} \\
\Gamma_1 & = \mathcal{A} \cap (\bar{B}_{q, 3}(0))^n \times \mathcal{A}.
\end{align*}
We observe that $\Gamma_1 \subset \Gamma_2 \subset \mathcal{X}.$
Proposition~\ref{thm:2} in Section~\ref{sec:result2} in combination with the fact that the right-hand sides of the $\dot{T}_i$ are well-defined, guarantees that $\mathcal{X}$ 
is forward invariant and can serve as the state space for 
$(x, T^{\text{tot}})$. Also the set $\Gamma_2$ is globally uniformly asymptotically stable relative to $\mathcal{X}$. 

On $\mathcal{X}$, the dynamics for $T_i$ is given by 
$$\dot{T}_i = \sum_{j\in\mathcal{N}_i}a_{ij}\left(\left(T_j - T_i\right) - R_iR_j^T\left(T_i -T_j\right)\right),$$
but on the set $\Gamma_2$ the dynamics for $T_i$ is given by 
$$\dot{T}_i = \sum_{j\in\mathcal{N}_i}a_{ij}\left(T_j - T_i\right).$$
The dynamics for $T^{\text{tot}}$ is given by 
$$\dot{T}^{\text{tot}} = -(L(\mathcal{G},A) \otimes I) T^{\text{tot}},$$
where $L(\mathcal{G},A)$ is the graph Laplacian matrix for 
a quasi-strongly connected graph. It is well known that the consensus set is 
exponentially stable for this dynamics. Thus, the set 
$\Gamma_1$ is globally asymptotically stable relative to $\Gamma_2$.
Now one can use Theorem 10 in \cite{reduction} in order to show that
$\Gamma_1$ is globally asymptotically stable relative to $\mathcal{X}.$
\hfill $\blacksquare$ \vspace{3mm}

\subsection{Numerical experiments}\label{sec:5.1}
\begin{figure}
\includegraphics[scale=0.160]{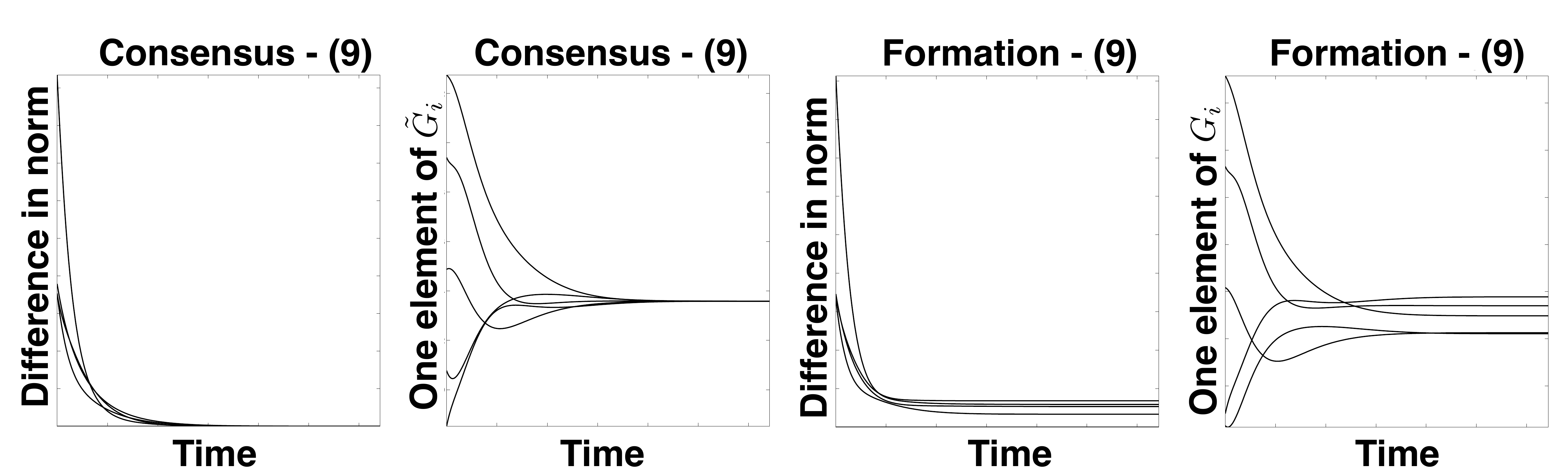}
\includegraphics[scale=0.160]{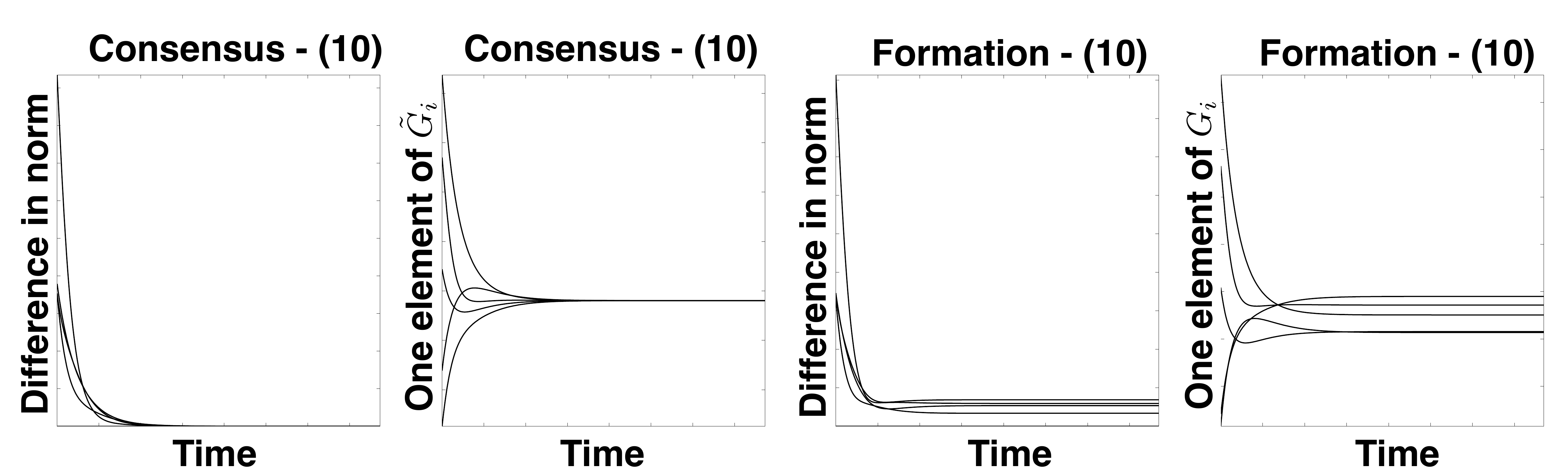}
\caption{These plots illustrate  the difference between reaching
consensus in the $\tilde{G}_i$ variables and reaching 
a formation in the $G_i$ variables. 
}
\label{fig:1}
\end{figure}

In order to illustrate the relation between consensus 
and formation the following example is considered. For a system of five agents,
in Figure~\ref{fig:1} the convergence of the $\tilde{G}_i$
variables to consensus and the 
convergence of the $G_i$ variables to a desired formation 
is shown. The adjacency matrix was chosen to that of a quasi-strongly connected graph with 
entries equal to $0$, $1$ or $2$.
The initial rotations are drawn from the uniform distribution 
over $B_{\pi/2}(I)$. 
Each initial translation vector is  drawn from the 
uniform distribution over the unit box in $\mathbb{R}^3$. The initial $R_i(0)$ rotations and initial $T_i(0)$ 
translations are the building blocks of the $G_i(0)$ transformations. The desired $G^*_i$ are constructed in 
the same manner as the $G_i(0)$, after which the $\tilde{G}_i(0)$ transformations
are constructed by $\tilde{G}_i(0) = {G}_i(0)G_i^{*-1}$.

For the same initial conditions, 
the four upper plots in Figure~\ref{fig:1} 
show the convergence when 
controller \eqref{eq:approach1:controller:1} is used, whereas
the four lower plots in Figure~\ref{fig:1} show the convergence when
 controller \eqref{eq:approach1:controller:2} is used. In each 
 of these four subplots, the first plot is showing the 
 difference $\|\tilde{G}_i(t) - \tilde{G}_1(t)\|_{\text{F}}$ 
 for all $i$; the second plot is 
 showing one of the elements of $\tilde{G}_i(t)$ for all $i$ as function of time, this element is the upper left one in the $\tilde{G}_i(t)$, \emph{i.e.}, it is an element of the rotation matrix;
 the third plot is showing the difference $\|{G}_i(t) - {G}_1(t)\|_{\text{F}}$ 
 for all $i$ as function of time; the fourth plot is 
 showing one of the elements of ${G}_i(t)$  as function of time for all $i$. This element is chosen as the upper left element in the ${G}_i(t)$.

The construction of the initial rotations in this example
 does not guarantee
that initial rotations are contained in the regions 
specified in Proposition~\ref{proposition:first_control_law:1}
and Proposition~\ref{proposition:first_control_law:2}, yet the convergence is obtained for both control laws. For 1000 simulations with five agents and random quasi-strongly connected 
topologies, where the initial rotations are drawn from the uniform distribution over 
$SO(3)$ and the translations are drawn from the uniform distribution over the unit cube
in $\mathbb{R}^3$,  the $\tilde{G}_i(t)$ transformations converged to consensus 909 respective 910 times for the two different control laws, \emph{i.e.}, a success rate of over 90 $\%$. If the initial 
rotations in $\tilde{G}_i$ were drawn from the uniform distribution over $B_{\pi/2}(I)$ the transformations converged to 
consensus 1000 respective 1000 times for the two different control laws, \emph{i.e.}, a success rate of 100~$\%$.

\begin{figure}
\includegraphics[scale=0.159]{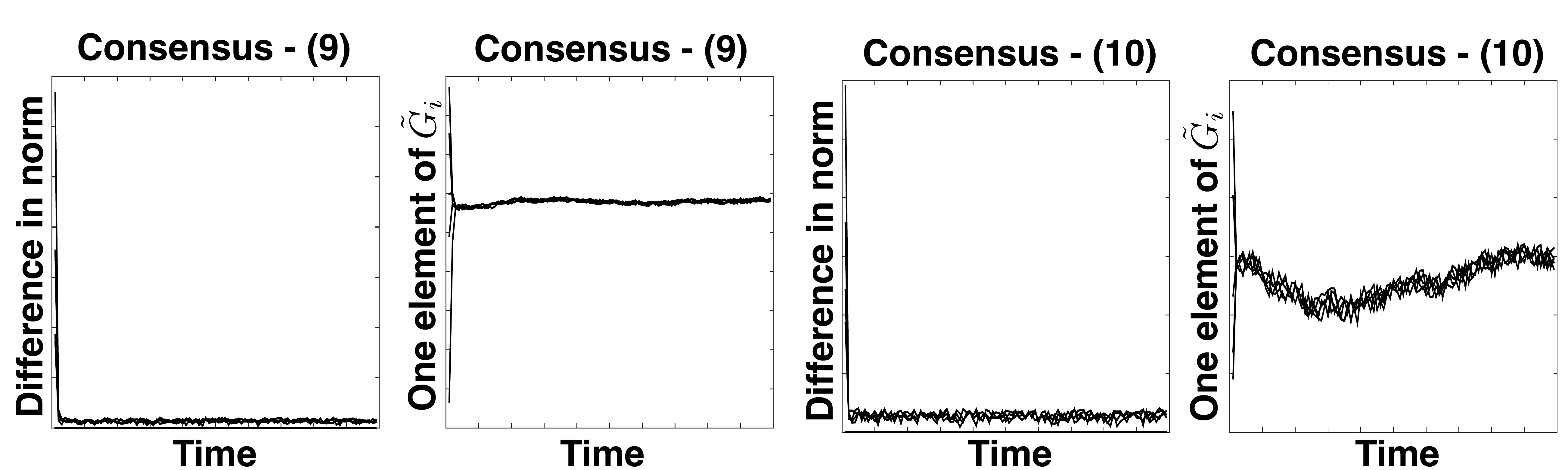}
\caption{These plots show the convergence to 
consensus under discrete sampling, additive absolute noise, and
switching between quasi-strongly connected graphs.
}
\label{fig:nisse1}
\end{figure}

Furthermore, numerical experiments were conducted 
when there was additive absolute, noise and when the transformations were measured discretely the 
topologies were switching.
The \enquote{noise} were random skew symmetric matrices, whose magnitudes
were equal to $0.1$.
Under these conditions the controllers \eqref{eq:approach1:controller:1} and \eqref{eq:approach1:controller:2}
were tested 
in 100 simulations where the graphs 
were switching between quasi-strongly connected topologies and the initial 
rotations in $\tilde{G}_i$ were drawn from the uniform distribution over $B_{\pi/2}(I)$. The number of agents was 5. The matrices converged to 
consensus in every simulation for both controller \eqref{eq:approach1:controller:1} and controller \eqref{eq:approach1:controller:2}. In the simulations the 
graphs switched with a frequency of $10$, which was the same as the
sampling frequency; the consensus is shown for one simulation in Figure~\ref{fig:nisse1}. 
The simulations show stronger results than those presented in 
Proposition~\ref{proposition:first_control_law:1} and 
Proposition~\ref{proposition:first_control_law:2}.

The input is constant between sample points. Thus, we can solve the system exactly between those points (it becomes a linear  time-invariant system). 
The solutions between the sampling points are not shown 
in the figure, instead there are straight lines 
connecting the solutions at the sample points. 

In all simulations the \enquote{random} graphs were created by constructing adjacency matrices in the following way: First an adjacency matrix for a tree graph 
was created and then a binary matrix was created where each element in the matrix was drawn from the uniform distribution over $\{0,1\}$. The final adjacency was then chosen as the sum of the adjacency matrix for the tree graph and the binary matrix.

\section{Results for the second control laws}\label{sec:result2}

\subsection{Rotations}
Here we address the controllers \eqref{chapter2:controller:1} and \eqref{chapter2:controller:2}. We start with 
\eqref{chapter2:controller:1}. 
The structure of controller \eqref{chapter2:controller:1} is well known from
the consensus problem in a system of agents with single integrator dynamics and states in $\mathbb{R}^m$~\cite{mesbahi2010graph}. 
The question is if this simple control law also works for rotations expressed in any of the local representations that we consider. The answer is yes. For all the convergence results provided in this section it is true that the state $x(t)$ converges to a fixed point, \emph{i.e.}, not a limit cycle.
\begin{prop}\label{thm:1}
Suppose $q < r$ and the graph $\mathcal{G}_{\sigma(t)}$ is uniformly strongly connected, then if controller \eqref{chapter2:controller:1} is used, 
$(\bar{B}_{q,3}(0))^n$ is strongly forward invariant, $0 \in \mathbb{R}^{mn}$ is uniformly stable
and $(\bar{B}_{q,3}(0))^n \cap \mathcal{A}$ is globally attractive relative to $(\bar{B}_{q,3}(0))^n$.
\end{prop}

In order to prove Proposition~\ref{thm:1},
we use the following proposition. 

\begin{prop}\label{prop:12}\label{prop:13}
Suppose controller \eqref{chapter2:controller:1} is used, $\mathcal{G}_{\sigma(t)}$ is uniformly strongly connected,
and $q < r$. 

Now, suppose there is a continuously differentiable function 
$$V: \mathbb{R}^3 \rightarrow \mathbb{R}$$
such that for any given  $k \in \{1, \ldots 2^{n-1}\}$ and 
$\bar{x} = [\bar{x}_1^T, \bar{x}_2^T, \ldots, \bar{x}_n^T]^T \in (\bar{B}_{q,3}(0))^n$
\begin{enumerate}
\item if 
$$i \in \textnormal{arg}\max_{j \in \mathcal{N}^{k}_i}(V(\bar{x}_j))$$ it holds that 
\begin{equation}\label{eq:t1}
\langle \nabla V(\bar{x}_i) , \sum_{j\in\mathcal{N}^{k}_i} a_{ij}({y}_{j}(\bar{x}) - {y}_{i}(\bar{x})) \rangle \leq 0,
\end{equation}
where $y_i$ is the local representation. \\
\item and equality holds for \eqref{eq:t1} if and only if $\bar{x}_i = \bar{x}_j$ for all $j \in \mathcal{N}^{k}_i$,
\end{enumerate} 
then $(\bar{B}_{q,3}(0))^n$ is strongly forward invariant for the dynamics of $x$ and $(\bar{B}_{q,3}(0))^n \cap \mathcal{A}$ is globally attractive
relative to $(\bar{B}_{q,3}(0))^n$. 
\end{prop}

The proof of Proposition~\ref{prop:12} is omitted here but follows,
up to small modifications, the 
procedure in the proof of Theorem~2.21. in \cite{thunberg2014distributed}. The essential difference between the two
is that besides the fact that in Theorem~2.21. in \cite{thunberg2014distributed} more general 
right-hand sides of the system dynamics are considered,
only one switching signal function is used for the system in that theorem,
whereas in in this work we assume individual switching signal functions 
for the agents.

\quad \emph{Proof of Proposition~\ref{thm:1}: }
We verify that \emph{(1)} and \emph{(2)} are satisfied in 
Proposition~\ref{prop:13} by choosing $V(\bar{x}_i) = \bar{x}_i^T\bar{x}_i$.
Let $i \in \textnormal{arg}\max_{j \in \mathcal{N}^{k}_i}(V(\bar{x}_j))$. Then
\begin{align*}
\langle \nabla V(\bar{x}_i),  \sum_{j\in\mathcal{N}^{k}_i} a_{ij}({y}_{j}(\bar{x}) - {y}_{i}(\bar{x})) \rangle & \leq \\
\sum_{j\in\mathcal{N}^{k}_i}(\|\bar{x}_j\|g(\|\bar{x}_j\|) - \|\bar{x}_i\|g(\|\bar{x}_i\|)) & \leq 0,
\end{align*}
where we have used the fact that
$g$ is strictly increasing. The last inequality is strict 
if and only if $\bar{x}_i = \bar{x}_j$ for all $j \in \mathcal{N}^{k}_i$.
\hfill $\blacksquare$ \vspace{3mm}

\begin{rem}
Instead of using \eqref{chapter2:controller:1}, one could use feedback linearization 
and construct the following control law for agent $i$,
$$\omega_i =  L_{y_i}^{-1}\sum_{j\in\mathcal{N}^{\sigma^i(t)}_i} a_{ij}({y}_{j} - {y}_{i}),$$
where $L_{y_i}$ is the Jacobian matrix for the representation $y_i$. If this feedback linearization
control law is used and the graph $\mathcal{G}_{\sigma(t)}$ is 
quasi-strongly connected,
the consensus set, restricted to any closed ball $(\bar{B}_{q,3}(0))^n$
where $q < r$, 
is globally uniformly asymptotically stable relative to $(\bar{B}_{q,3})^n$. However, for many representations such as the Rodrigues Parameters, the Jacobian matrix $L_{y_i}$
is close to singular as $y_i$ is close to the boundary of $\bar{B}_{q,3}(0)$. Furthermore, the expression is nonlinear in the $y_i$. This might make this type of
control law more sensitive to measurement errors than \eqref{chapter2:controller:1}. 
\end{rem}

Now we continue with the study of (\ref{chapter2:controller:2}) where only local representations of the relative
rotations are available. Under stronger assumptions on the
initial rotations of the agents at time $t_0$ and weaker assumptions on
the graph $\mathcal{G}_{\sigma(t)}$, the following proposition ensures
uniform asymptotic convergence to the consensus set.

\begin{prop} \label{thm:2}
Suppose 
$q < r/2$ and the controller \eqref{chapter2:controller:2} is used, then 
$(\bar{B}_{q,3}(0))^n$ is strongly forward invariant and $(\bar{B}_{q,3}(0))^n \cap \mathcal{A}$ is globally uniformly asymptotically stable relative to 
$(\bar{B}_{q,3}(0))^n$ if and only if $\mathcal{G}_{\sigma(t)}$ is uniformly quasi-strongly connected.
\end{prop}

\begin{rem}
In Proposition \ref{thm:2}, since only information that is independent of $\mathcal{F}_W$ is used in 
\eqref{chapter2:controller:2}, the assumption that the rotations initially are contained in $\bar{B}_q(I)$ can be relaxed.
As long as there is a $Q \in SO(3)$ such that all the rotations 
are contained in $(\bar{B}_{q}(Q))^n$ initially, 
the rotations will reach consensus asymptotically and uniformly with respect to time.
\end{rem}
In order to prove Proposition \ref{thm:2}, we first provide 
a theorem, which gives some geometric insight. Then we provide a
Proposition, which guarantees asymptotic stability of the 
consensus set.

\begin{thm}\label{thm:3}
Suppose that the control law \eqref{chapter2:controller:2} is used and $x \in (B_{q, 3}(0))^n$, where
$q < r/2$. Let
$z_i = \tan(\theta_i/2)u_i$ and $z = [z_1^T, \ldots, z_n^T]^T$.
Then
\begin{align*}
\dot{z}_1  = & \sum_{j \in \mathcal{N}^{\sigma^1(t)}_1}a_{1j}h_{1j}(z_1,z_j)(z_{j} - z_1), \\
\vdots & \\
\dot{z}_n  = & \sum_{j \in \mathcal{N}^{\sigma^n(t)}_n}a_{nj}h_{nj}(z_n,z_j)(z_{j} - z_n),
\end{align*}
where $h_{ij}(z_i,z_j) \geq 0$ and $h_{ij}(z_i,z_j) > 0$ if $z_{j} \neq z_i$. 
\end{thm}

\begin{rem}
The $h_{ij}$ functions in Theorem~\ref{thm:3} depends on the parameterization $y$.
\end{rem}

A proof of Theorem~\ref{thm:3} (up to small modifications
due to the assumptions on the switching signal functions)
can be found in~\cite{thunberg2014distributed}. It is based 
on the results in \cite{afsari2011riemannian,Hartley2011}.
Theorem \ref{thm:3} states that, after 
a change of coordinates to the Rodrigues Parameters,
the system satisfies the well known convexity 
assumption that the right-hand side of 
each agent's dynamics is inward-pointing~\cite{afsari2011riemannian} relative
to the convex hull of its neighbors' positions. There 
are many publications addressing this type of dynamics, \emph{e.g.}, \cite{moreau2005stability,shi2009,lin2007state}.

\begin{prop}\label{prop:18}
Suppose control law \eqref{chapter2:controller:2} is used,  $q < r/2$, and $(\bar{B}_{q,3}(0))^n$ is strongly forward invariant for the dynamics of $x$. 

Suppose there is a continuously differentiable function 
$$W: \mathbb{R}^3 \times \mathbb{R}^3 \rightarrow \mathbb{R}^+,$$
such that for any given  $k,l \in \{1, \ldots 2^{n-1}\}$ and 
$\bar{x} = [\bar{x}_1^T, \bar{x}_2^T, \ldots, \bar{x}_n^T]^T \in (\bar{B}_{q,3}(0))^n$, 
\begin{enumerate}
\item if 
$$(i,j) \in \textnormal{arg}\max_{k' \in \mathcal{N}^{k}_i, l' \in \mathcal{N}^{l}_j}(W(\bar{x}_i,\bar{x}_j))$$ it holds that 
\begin{align*}
\langle \nabla W(\bar{x}_i,\bar{x}_j), [\left(\sum_{k'\in\mathcal{N}^{k}_i} a_{ij}({y}_{k'}(\bar{x}) - {y}_{i}(\bar{x}))\right)^T, & \\
\left(\sum_{l'\in\mathcal{N}^{l}_j} a_{ij}({y}_{l'}(\bar{x}) - {y}_{j}(\bar{x}))\right)^T] \rangle & \leq 0 
\end{align*}
\item and equality holds if and only if $\bar{x}_i = \bar{x}_{k'}$ for all $k' \in \mathcal{N}^{k}_i$
and $\bar{x}_j = \bar{x}_{l'}$ for all $l' \in \mathcal{N}^{l}_j$,
\end{enumerate}
then $(\bar{B}_{q,3}(0))^n \cap \mathcal{A}$ is globally uniformly asymptotically stable
relative to $(\bar{B}_{q,3}(0))^n$ if and only if
$\mathcal{G}_{\sigma(t)}$ is uniformly quasi-strongly connected.
\end{prop}

The proof of Proposition~\ref{prop:18}, is omitted here, but follows,
up to small modifications, the 
procedure in the proof of Theorem~2.22. in \cite{thunberg2014distributed}.

\vspace{2mm}
\quad \emph{Proof of Proposition~\ref{thm:2}}:
Let us define the functions
\begin{align*}
V(x_i) = & \: \: \: \: x_i ^Tx_i \text{ and } \\
W(x_i,x_j) = & (z_j(x_j) - z_i(x_i))^T(z_j(x_j) - z_i(x_i)),
\end{align*}
Using Theorem~\ref{thm:3} and the function $V$ together with Proposition~\ref{prop:12}, along the lines of the 
proof of Proposition~\ref{thm:1}, one can show that 
$(\bar{B}_{q,3}(0))^n$ is strongly forward invariant 
and if $\mathcal{G}_{\sigma(t)}$ is uniformly strongly connected, $(\bar{B}_{q,3}(0))^n \cap \mathcal{A}$ is 
globally attractive relative to $(\bar{B}_{q,3}(0))^n$. 

Now, since $(\bar{B}_{q,3}(0))^n$ is strongly forward invariant, one can use Theorem ~\ref{thm:3} in order to  show that $W$ satisfies the criteria in Proposition~\ref{prop:18}. The mapping 
$$x_i \mapsto z_i$$
is a diffeomorphism on $(\bar{B}_{r,3}(0))^n$.
The set $(\bar{B}_{q,3}(0))^n \cap \mathcal{A}$ is globally 
uniformly asymptotically stable relative to $(\bar{B}_{q,3}(0))^n$.
\hfill $\blacksquare$ \vspace{3mm}

\begin{rem}
We can generalize the results in 
Proposition~\ref{thm:1} and Proposition~\ref{thm:2}.
Up until now we have assumed that we first 
fix a representation $y_i$, $y_{ij}$ and then we use the control laws
\eqref{chapter2:controller:1} and \eqref{chapter2:controller:2} for this representation. 
Instead, at each switching time $\tau_k^i$ we can allow
the representation to switch also. 
\end{rem}

The following proposition addresses a special case when the 
the rate of convergence is exponential.

\begin{prop}\label{proposition:exponential}
Suppose $\mathcal{G}_{\sigma(t)}$ fulfills the following. At each 
time $t$ and for each pair $(i,j)$, the edge $(i,j) \in \mathcal{E}_{\sigma(t)}$ or the edge $(j,i) \in \mathcal{E}_{\sigma(t)}$.  
Suppose controller \eqref{chapter2:controller:2} is used
and $g(\theta_i) \geq k\theta_i$ for some $k > 0$.
For $q < r/2$, the set 
$\{(R_1, \ldots, R_n) \in (\bar{B}_{q}(I))^n: R_1 = \hdots = R_n\}$ is globally exponentially stable
relative to $(\bar{B}_{q}(I))^n$ for the closed loop dynamics of $(R_1, R_2, \ldots, R_n)$ with respect to the Riemannian metric on $SO(3)$.
\end{prop}

\begin{rem}
Using the results in \cite{Chen:2014hj} it can be obtained that for the 
Axis-Angle Representation the convergence is exponential also for general
uniformly quasi-strongly connected graphs, {i.e.}, not only the 
restricted class of graphs considered Proposition~\ref{proposition:exponential}.
\end{rem}

In Proposition~\ref{proposition:exponential} \emph{(1)}, since we have assumed that $g$ is analytic, the condition $g(\theta_i) \geq k\theta_i$ can equivalently be formulated as 
$g(\theta_i) = \mathcal{O}(\theta_i)$ as $\theta_i \rightarrow 0$. All the local representations 
previously addressed fulfill this assumption, \emph{e.g.}, the 
Axis-Angle Representation, the Rodrigues Parameters and  the Unit Quaternions.

Before we prove Proposition~\ref{proposition:exponential} we 
formulate the following lemma.

\begin{lem}\label{lemma:support}
Suppose $x \in (\bar{B}_{q,3}(0))^n$ where $q < r/2$. If 
$$(i,j) \in \textnormal{arg}\max_{(k,l) \in \mathcal{V} \times \mathcal{V}}\|x_{kl}\|,$$
then
$$x_{ij}^Ty_{ik} \geq 0 \quad \text{ for all } \quad k.$$ 
\end{lem}

\vspace{2mm}
The proof of Lemma~\ref{lemma:support} follows more or less as a consequence  of Theorem~\ref{thm:3} and is omitted here.

\vspace{2mm}
\quad \emph{Proof of Proposition~\ref{proposition:exponential}}:
We already know from Proposition~\ref{thm:2} that the set 
$(\bar{B}_{q,3}(0))^n$ is strongly forward invariant and
$(\bar{B}_{q,3}(0))^n \cap \mathcal{A}$ is globally uniformly 
asymptotically stable relative to $(\bar{B}_{q,3}(0))^n$. What is 
left to  prove is that for the special structure of the graph considered,
the rate of convergence is exponential relative to $(B_{q}(I))^n$ when the Riemannian metric is used. 

Let us define $$\alpha = \min_{(k,l) \in \mathcal{V} \times \mathcal{V}}a_{kl},$$
and $$V(x) = \max_{(k,l) \in \mathcal{V} \times \mathcal{V}}x_{kl}^Tx_{kl}.$$
At time $t$ let $(i,j)$ be such that $V(x(t)) = x^T_{ij}(t)x_{ij}(t)$.
$$x_{ij}^T\left(\omega_j - \omega_i\right) \leq -\alpha k V(x(t)),$$
where the last inequality is due to Lemma~\ref{lemma:support} and the assumption on the 
graph $\mathcal{G}_{\sigma(t)}$.
Now one can show that 
$$D^+V(x(t)) \leq -\alpha k V(x(t)).$$
By using the Comparison Lemma, one can show that $V$ converges to 
zero with exponential rate of convergence. 
\hfill $\blacksquare$ \vspace{3mm}

\subsection{Illustrative example}\label{sec:6.1}
\begin{figure}
\includegraphics[scale=0.23]{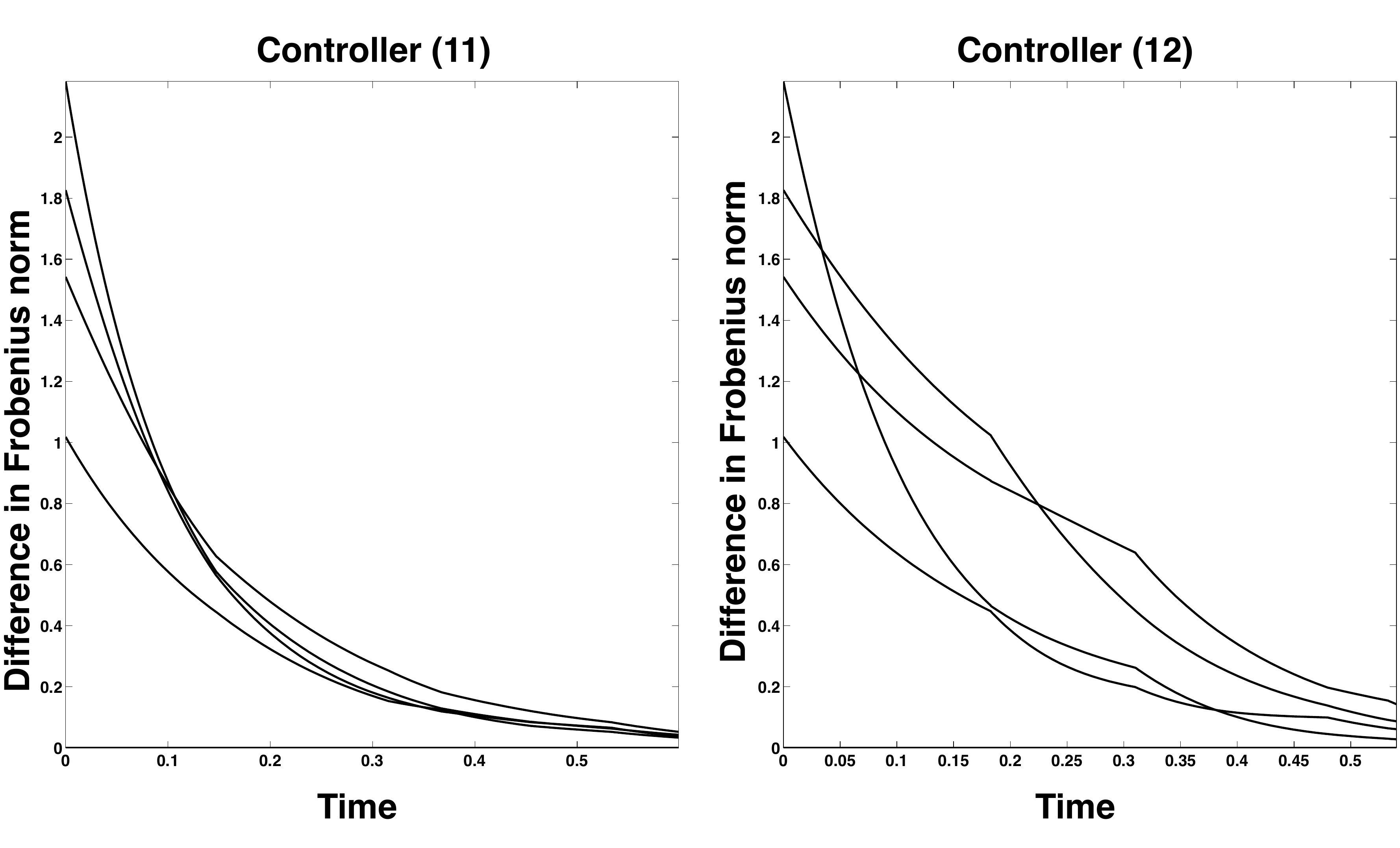}
\caption{Convergence to consensus. For the same initial 
rotations, controller \eqref{chapter2:controller:1} and 
controller \eqref{chapter2:controller:2} are used. The left
plot shows the errors (in terms of Frobenius norm) between the
first rotation and the other rotations  as a function of time
when controller \eqref{chapter2:controller:1} is used;
in the right plot the same type of errors is shown when controller \eqref{chapter2:controller:2} is used.
}
\label{figure:2}
\end{figure}

In order to illustrate the convergence of the 
rotations to the consensus set, an 
illustrative example is constructed where the representation $(R - R^T)^{\vee}$ is chosen 
both for control law \eqref{chapter2:controller:1}
and \eqref{chapter2:controller:2}. The number 
of agents is $5$ and the graph the 
graphs were constructed in the same manner as
in Section~\ref{sec:5.1}.  The initial 
rotations are drawn from the uniform distribution over $B_{\pi/2}(I)$. The convergence to consensus is shown in Figure~\ref{figure:2}.

\subsection{Translations}
Here we address the controllers \eqref{chapter2:controller:1:v} and \eqref{chapter2:controller:2:v}

Controller \eqref{chapter2:controller:1:v}, despite 
its appealing structure does in general not guarantee
consensus in the translations. In order to see this,
we consider the following example. 

$$\dot{T}_i = \sum_{j \in \mathcal{N}^{\sigma^i(t)}_i}R_i^T(T_j -T_i).$$
Suppose 
\begin{align*}
& \begin{bmatrix}
1 & 0 & 0 \\
0 & 1 & 0 \\
0 & 0 & 0
\end{bmatrix}T_i(0) = T_i(0) \text{ and } \\
R_i(t) = & 
\begin{bmatrix}
-1 & 0 & 0 \\
0 & -1 & 0 \\
0 & 0 & 1
\end{bmatrix} \quad \text{ for all } i,t.
\end{align*}
Then, for this particular choice of rotations and initial conditions
for the $T_i$,
$$\dot{T}_i = \sum_{j \in \mathcal{N}^{\sigma^i(t)}_i}(T_i -T_j).$$
Thus,
$$\dot{T}^{\text{tot}} = (L(\mathcal{G}_{\sigma(t)}, \mathcal{A}) \otimes I)T^{\text{tot}},$$
which is unstable. 

One partial result for controller \eqref{chapter2:controller:1:v} is the following one.
By a change of coordinates one can prove that,
if all the rotations of the agents are the same
and constant and the translations are contained in 
the linear subspace spanned by the rotational axis, 
the translations converge asymptotically to consensus.

Controller \eqref{chapter2:controller:2:v} delivers
a much stronger result. 

\begin{prop}
Suppose controller \eqref{chapter2:controller:2:v} is used
and the graph $\mathcal{G}_{\sigma(t)}$
is uniformly quasi-strongly connected. The set $\mathcal{A}$ is globally asymptotically stable stable.
\end{prop}
The proof of the proposition is based on the fact that the closed loop dynamics
is
given by
$\dot{T}^{\text{tot}} = -L(\mathcal{G}_{\sigma(t)},A)T^{\text{tot}}.$

\section{Control on the dynamic level for rigid bodies in space}\label{sec:dyn}
In this section we construct control laws on the dynamic level 
for the case of rigid bodies in space. 
The dynamical equations for agent $i$ are given by
\begin{equation}\label{eq:s1}
\begin{cases}
\dot{R}_i  = R_i\widehat{\omega}_i, \\
\dot{\omega}_i  = J_i^{-1}\left(-\widehat{\omega}_iJ_i\omega_i + \boldsymbol{\tau}_i\right), \\
\dot{T}_i  = R_iv_i, \\
\dot{v}_i  = \frac{\boldsymbol{f}_i}{m_i} - \widehat{\omega}_iv_i,
\end{cases}
\end{equation}
where $J_i$ is the inertia matrix and $\boldsymbol{\tau}_i$ is the 
control torque; $\omega_i$ is a state variable.
In the formation problem the goal is to reach consensus in the 
$\tilde{(\cdot)}$-variables, and the dynamical equations 
for those variables are
\begin{equation}\label{eq:s2}
\begin{cases}
\dot{\tilde{R}}_i &  = \tilde{R}_i\widehat{\tilde{\omega}}_i, \\
\dot{\tilde{\omega}}_i & =  R_i^{*}J_i^{-1}\left(-\left(R_i^{*T}\widehat{\tilde{\omega}}_i R_i^{*}\right)J_iR_i^{*T}\tilde{\omega}_i + \boldsymbol{\tau}_i\right), \\
\dot{\tilde{T}}_i & = \tilde{R}_i\tilde{v}_i \\
\dot{\tilde{v}}_i & = R_i^{*}\left (\frac{\boldsymbol{f}_i}{m_i} - R_i^{*T} \widehat{\tilde{\omega}}_i\left (\widehat{\tilde{\omega}}_iT_i^* + \tilde{v}_i\right )\right) - \left (\dot{\tilde{\omega}}_i\right )^{\wedge} T_i^*.
\end{cases}
\end{equation} 

In this section, we strengthen the assumptions on 
$\mathcal{G}_{\sigma(t)}$ by assuming it is time-invariant. Thus, we 
 denote the time-invariant (also referred to as constant or fixed) graph by $\mathcal{G}$. 
The reason for choosing time-invariant graphs
is that we are now considering a second order system, and the methods we use here
are based on backstepping. In order to show stability, we introduce auxiliary error variables,
and in the case of a switching graph, these variables suffer from discontinuities. 
One way to avoid this problem is to replace the
discontinuities with continuous in time transitions. This is however
not something we do here.

\subsection{Rotations}
Only the consensus problem and the first set
of equations, \eqref{eq:s1}, will be considered here. When performing 
formation control, the presented control laws below, \eqref{chapter2:torque:controller:3} and \eqref{chapter2:torque:controller:4}, are modified slightly. 
In both control laws, all the variables should be replaced
by $\tilde{(\cdot)}$-variables, \emph{i.e.}, $x_i$ should be $\tilde{x}_i$ instead, $\tilde{\omega}_i$ should be ${\omega}_i$ instead and so on. The expression \enquote{$J_i$(} is replaced by \enquote{$J_iR_i^T$(},
and the expression \enquote{$\widehat{\omega}_iJ_i\omega_i$}
is replaced by \enquote{$\left(R_i^{*T}\widehat{\tilde{\omega}}_i R_i^{*}\right)J_iR_i^{*T}\tilde{\omega}_i$}.

Based on the two kinematic control laws \eqref{chapter2:controller:1} and \eqref{chapter2:controller:2}, we now propose two torque control laws
for each agent $i$, 
where the first one is based on absolute rotations and the second one is based
on relative rotations. The control laws are
\begin{align}
\label{chapter2:torque:controller:3}
\boldsymbol{\tau}_i & = J_i(-x_i +  \sum_{j \in \mathcal{N}_i}a_{ij}(L_{x_j}\omega_j - L_{x_i}\omega_i - \bar{\omega}_i)) + \widehat{\omega}_iJ_i\omega_i, \\
\label{chapter2:torque:controller:4}
\boldsymbol{\tau}_i & = J_i( -k_i\bar{\omega}_i' + \sum_{j \in \mathcal{N}_i}a_{ij}L_{-y_{ij}}\omega_{ij}) + \widehat{\omega}_iJ_i\omega_i.
\end{align}
The parameter $k_i$ is a positive gain. The error variables $\bar{\omega}_i$ and $\bar{\omega}_i'$ are by
follows
\begin{align*}
\bar{\omega}_i & = \omega_i  - \sum_{j \in \mathcal{N}_i}a_{ij}(x_j - x_i), \\
\bar{\omega}_i' & = \omega_i  - \sum_{j \in \mathcal{N}_i}a_{ij}y_{ij}.
\end{align*}
The matrix $L_{y_{ij}}$ is the Jacobian matrix for 
$y_{ij}$, \emph{i.e.},
\begin{align*}
\dot{y}_{ij} & = L_{-y_{ij}}\omega_{ij},
\end{align*}
and
$$\omega_{ij} = R_{ij}\omega_j - \omega_i$$ is the relative
angular velocity between agent $i$ and agent $j$.  
In the following, the notation $(x_i,\bar{\omega}_i') = [x_i^T,\bar{\omega}_i'^T].$  
We collect all the $x_i$ and $\bar{\omega}_i$ into $(x, \bar{\omega}) \in (B_{\pi,3})^n \times (\mathbb{R}^3)^n$
and all the $x_i$ and $\bar{\omega}_i'$ into $(x, \bar{\omega}') \in (B_{\pi,3})^n \times (\mathbb{R}^3)^n$.
Now, given $i \in \mathcal{V}$, the
right-hand side for $(x_i, \bar{\omega}_i)^T$
when the torque control law \eqref{chapter2:torque:controller:3} is used is
\begin{align*}
\dot{x}_i & = L_{x_i}\sum_{j \in \mathcal{N}_i}a_{ij}(x_j - x_i) + L_{x_i}\bar{\omega}_i, \\
\dot{\bar{\omega}}_i & = -x_i - \sum_{j \in \mathcal{N}_i}a_{ij}\bar{\omega}_i,
\end{align*}
whereas the closed loop system for $(x_i, \bar{\omega}'_i)^T$ when the torque 
control law \eqref{chapter2:torque:controller:4} is used, is
\begin{align*}
\dot{x}_i & = L_{x_i}\sum_{j \in \mathcal{N}_i}a_{ij}y_{ij} + L_{x_i}\bar{\omega}_i', \quad \quad \: \: \: \\
\dot{\bar{\omega}}_i' & = -k_i\bar{\omega}_i'.
\end{align*}

We note that in \eqref{chapter2:torque:controller:3}, each agent $i$ needs to know, not only the absolute
rotations of its neighbors, but also the angular velocities of its neighbors. This requirement 
is fair, in the sense that in order to obtain the absolute rotations of the neighbors,
communication is in general necessary. In this case the angular velocities can also 
be transmitted. In \eqref{chapter2:torque:controller:4}, we see that each agent $i$ needs to know the relative 
rotations, relative velocities to its neighbors and the angular velocity of itself. The assumption that agent $i$  
knows its own angular velocity is quite strong the sense that this velocity is not to be regarded as relative information.
However, in practice the angular velocity is possible to measure without the knowledge of the 
global frame $\mathcal{F}_W$. Thus, the angular velocity is local information.

\begin{prop}\label{thm:4}
Suppose $\mathcal{G}$ is strongly connected.
If $$\max_{i \in \mathcal{V}} \: x_i^T(t_0)x_i(t_0) + \bar{\omega}_i(t_0)^T\bar{\omega}_i(t_0) \leq q < \pi,$$
{i.e.}, $(x_i(t_0), \bar{\omega}_i(t_0))^T \in \bar{B}_{q,6}$ for all~$i$ and some $q < \pi$, then if controller \eqref{chapter2:torque:controller:3} is used,
$\bar{B}_{q,6}$ is invariant for $(x(t), \bar{\omega}(t))$ and $x(t) \rightarrow \mathcal{A}$ and $\omega_i(t) \rightarrow 0$ for all
$i$ as $t \rightarrow \infty$. 
\end{prop}

\quad \emph{Proof}:
In the multi-agent system at hand we have $n$ agents, where each agent $i$ has the state
$(x_i, \bar{\omega}_i)^T$. We first show the invariance of the ball $\bar{B}_{q, 6}$.
$$V((x_i, \bar{\omega}_i)^T) = \frac{1}{2}\left (x_i^Tx_i + \bar{\omega}_i^T\bar{\omega}_i \right ).$$
We see that 
\begin{align*}
& \frac{d}{dt}V((x_i, \bar{\omega}_i)^T) \\ 
 = & \sum_{j \in \mathcal{N}_i}a_{ij}(x_i, \bar{\omega}_i)(x_j - x_i, -\bar{\omega}_i)^T \\
 = & \sum_{j \in \mathcal{N}_i}a_{ij}((x_i, \bar{\omega}_i)(x_j, 0)^T - (x_i, \bar{\omega}_i)(x_i, \bar{\omega}_i)^T) \\
 \leq & \sum_{j \in \mathcal{N}_i}a_{ij}x_i^T(x_j - x_i).
\end{align*}
Thus, $$D^+f_{V,6}((x(t), \bar{\omega}(t))^T) \leq 0.$$
Now, by using the Comparison Lemma one can show the invariance. 

In order to show the convergence, we define the following function
$$\bar{\gamma}(x, \bar{\omega}) = \sum_{i= 1}^n\xi_i (x_i^Tx_i + \bar{\omega}_i^T\bar{\omega}_i),$$
where $\xi = (\xi_1, \ldots, \xi_n)^T$ is the positive vector chosen such that (the symmetrical part of) $\text{diag}(\xi)L(\mathcal{G},A)$ is positive semi-definite.
We have that 
$$\dot{\bar{\gamma}} =  -{x}^T({L}' \otimes {I}_3){x} - \sum_{i= 1}^n\xi_i \sum_{j \in \mathcal{N}_i}a_{ij}\bar{\omega}_i^T\bar{\omega}_i.$$
By LaSalle's theorem, $(x(t), \bar{\omega}(t))^T$ will converge to the largest invariant set contained
in $$\{(x, \bar{\omega})^T: \dot{\bar{\gamma}}((x, \bar{\omega})^T)~= 0\}$$ as the time goes to infinity. This largest invariant set is contained in the set $\{(x, \bar{\omega})^T: x \in \mathcal{A}, \bar{\omega} = 0\}$.
\hfill $\blacksquare$ \vspace{3mm}

\begin{rem}
In the proof of Proposition ~\ref{thm:4}, if we look at the dynamics of $(x, \bar{\omega})$, we see that the largest invariant set contained in $\{(x, \bar{\omega})^T: \dot{\bar{\gamma}}((x, \bar{\omega})^T)~= 0\}$ is actually the point $0$. Hence, the system will reach consensus in the 
point $x = 0$.
\end{rem}

Now let us turn to control law~\eqref{chapter2:torque:controller:4}.
\begin{prop}\label{thm:5}
Suppose $\mathcal{G}$ is quasi-strongly connected.
For any positive $r_1$ and $r_2$ such that $r_1 < r_2 < r/2$ and $q >0$,
 there is a $k >0$ such that if $k_i \geq k$ and
$(x_i(t_0), \bar{\omega}_i'(t_0))^T \in \bar{B}_{r_1,3} \times \bar{B}_{q,3}$ for all $i$, then if controller \eqref{chapter2:torque:controller:4} is used 
it holds that
$(x_i(t), \bar{\omega}_i'(t))^T \in \bar{B}_{r_2,3} \times \bar{B}_{q,3}$ for all $i, t \geq t_0$ and
$$(x(t), \bar{\omega}'(t))^T \rightarrow (\bar{B}_{r_2,3})^n \cap \mathcal{A} \times \{0\} \quad \text{ as } \quad  t \rightarrow \infty.$$
Furthermore, $(\bar{B}_{r_2,3})^n \cap \mathcal{A} \times \{0\}$ is globally asymptotically stable relative to the largest
invariant set contained in $(\bar{B}_{r_2,3})^n \times (\bar{B}_{q,3})^n$ for the dynamics of $(x(t), \bar{\omega}'(t))^T$.
\end{prop}

\quad \emph{Proof of Proposition~\ref{thm:5}}:
Let us define
$$\mathcal{D}^* \subset \mathcal{D} = (\bar{B}_{r_2,3})^n \times (\bar{B}_{q,3})^n,$$
as the largest invariant set 
contained in $\mathcal{D}$. The set $\mathcal{D}^*$ is compact
and implicitly a function of $k$ (or the $k_i$).

Now we show that for a proper choice of the constant $k$, it holds that 
$$(\bar{B}_{r_1,3})^n \times (\bar{B}_{q,3})^n \subset \mathcal{D}^*.$$
We assume without loss of generality that $t_0 = 0$, and note that
$$\|\bar{\omega}_i'(t)\| =\|\bar{\omega}_i'(0)\|\exp(-k_it) \leq q\exp(-k_it) \leq q\exp(-kt).$$
We choose 
$$V(x_i(t)) = x_i^T(t)x_i(t).$$ 
By using Lemma~\ref{lemma:support}, it is possible to show that there exists an interval $[0, t_1)$ on which it holds that
$$D^+(\max_iV(x(t))) \leq qr_2\exp(-kt).$$
By using the Comparison Principle, it follows that 
$$\max_iV(x_i(t)) \leq \max_iV(x_i(0)) + qr_2\frac{(1- \exp(-kt))}{k}$$
on $[0, t_1)$.
Now if we choose $k \geq qr_2/(r_2 - r_1)$ we see that $\max_iV(x_i(t)) \leq r_2$
for $t \geq 0$, and we can choose $t_1 = \infty$. 

In order to show the desired convergence we use 
Theorem~10 in \cite{reduction}, where $\mathcal{X} = \mathcal{D}^*$,
$\Gamma_2 = \mathcal{D}^* \cap ((\bar{B}_{r_2,3})^n \times \{0\})$ and $\Gamma_1 = \mathcal{D}^* \cap (\mathcal{A} \times \{0\})$. 
\hfill $\blacksquare$ \vspace{3mm}

\subsection{Translations}
Also in this section only the consensus problem and the first set
of equations, \eqref{eq:s1}, are considered here.
We introducing a generalized version of the control law
\eqref{chapter2:controller:2:v}. 
When the formation control 
problem is considered the all variables are 
replaced by $\tilde{(\cdot)}$-variables. Furthermore 
the expression \enquote{$m_i($} is replaced by 
\enquote{$m_iR_i^{*T}($}, and the expression \enquote{$\widehat{\omega}_iv_i$}
is replaced by \enquote{$\widehat{\tilde{\omega}}_i(\widehat{\tilde{\omega}}_iT_i^* + \tilde{v}_i)) + (\dot{\tilde{\omega}}_i)^{\wedge} T_i^*$}.

The proposed consensus controller
is
\begin{align*}
\boldsymbol{f}_i  = & m_i(- k_i\bar{v}_i + \sum_{j\in\mathcal{N}_i} a_{ij}R_i^T(R_jv_j - R_iv_i) \\
& - \sum_{j\in\mathcal{N}_i} a_{ij}\widehat{\omega}_i R_i^T(T_j - T_i) + \widehat{\omega}_iv_i),
\end{align*}
where 
$$\bar{v}_i = v_i - \sum_{j\in\mathcal{N}_i} a_{ij}{T}_{ij}. \hspace{5mm}$$
The closed loop dynamics is
\begin{align*}
\hspace{10mm} \dot{T}_i & = \sum_{j\in\mathcal{N}_i} a_{ij}(T_j - T_i) + R_i(t)\bar{v}_i \\
\dot{\bar{v}}_i & -k_i\bar{v}_i.
\end{align*}
By treating the time as a variable $z$, we get the following system
\begin{align*}
\dot{z} & = 1 \\
\hspace{10mm} \dot{T}_i & = \sum_{j\in\mathcal{N}_i} a_{ij}(T_j - T_i) + R_i(z)\bar{v}_i \\
\dot{\bar{v}}_i & -k_i\bar{v}_i.
\end{align*}
Let the state of the entire system be $(z,T^{\text{tot}},\bar{v}^{\text{tot}})$, 
where $v^{\text{tot}} = [\bar{v}_1^T(t), \bar{v}_2^T(t), \ldots, \bar{v}_n^T(t)]^T \in \mathbb{R}^{3n}$.

\begin{prop}
Suppose that $R_i(z)$ is well behaved, in the sense that the right-hand side of the dynamics for $(z,T^{\text{tot}},{v}^{\text{tot}})$
is locally Lipschitz, then
the set $\mathbb{R} \times \mathcal{A} \times 0$ is globally asymptotically stable for the system.
\end{prop}

\quad \emph{Proof}:
Let the state space be $\mathcal{X} = \mathbb{R} \times \mathbb{R}^{3n} \times \mathbb{R}^{3n}$. We define the two closed subsets $\Gamma_1 \subset \Gamma_2$ of $\mathcal{X}$ as follows
\begin{align*}
\Gamma_1 & = \mathbb{R} \times \mathcal{A} \times 0 \\
\Gamma_2 & = \mathbb{R} \times \mathbb{R}^{3n} \times 0 \\
\end{align*}
It is easy to show that $\Gamma_2$ is globally asymptotically stable 
relative to $\mathcal{X}$ and $\Gamma_1$ is globally asymptotically 
stable relative to $\Gamma_2$.
Now the desired result follows from Theorem 10 in \cite{reduction}.
\hfill $\blacksquare$ \vspace{3mm}

\subsection{Illustrative examples}
\begin{figure}
\includegraphics[scale=0.160]{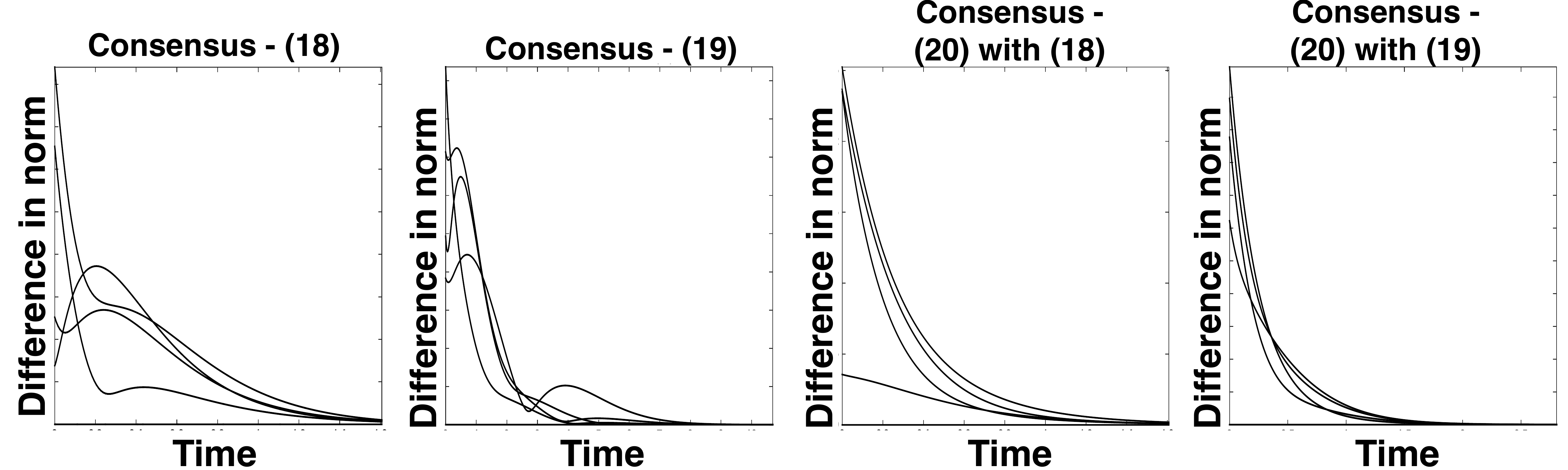}
\caption{Convergence to consensus. 
}
\label{figure:5}
\end{figure}

In Figure \ref{figure:5} the
convergence to consensus is shown when controllers
(18), (19) and (20) are used. In the simulations, five agents were considered and a random quasi-strongly connected graph was used.
The convergence 
to consensus is shown for the rotations, left plots,
and the translation, right plots, when controller
(18) was used together with controller (20) and controller
(19) was used together with controller (20).
The left plots shows the Euclidean distance between
$x_i(t)$ and $x_1(t)$ 
for $i = 2, \ldots, 5$, and the right plots 
show the Euclidean distance between
$T_i$ and $T_1$ as a function of time 
for $i = 2, \ldots, 5$. 

In controller (19) as well as controller (20)
the $k_i$ were chosen to $3$ for all $i$. 
The adjacency matrix was chosen to that of a quasi-strongly connected graph with 
entries equal to $0$, $1$ or $2$.
In controller (19) the representation
$(R - R^T)^{\vee}$ was used as the local representation for the $y_{ij}$.

\section*{Conclusions}
This work has considered the consenus and formation problems
on $SE(3)$ for multi-agent systems with switching 
interaction topologies. By a change of coordinates
it was shown that the consensus problem
can be seen as equivalent to the formation problem.
Any control law designed for the consensus problem can, after 
change of coordinates, be used for the formation problem.
New kinematic control laws 
have been presented as well as new convergence results. It has been shown that the same type of
control laws can be used for many popular local representations of $SO(3)$ such as the Modified Rodrigues Parameters and the Axis-Angle Representation. It has been shown that some 
of the control laws guarantee almost global convergence. For non-switching topologies, the 
kinematic control laws have been extended to torque
and force control laws for rigid bodies in space. 
The proposed control approaches have been justified 
by numerical simulations. 

\bibliographystyle{unsrt}     
\bibliography{Thesis}

\end{document}